\documentclass[12pt]{article}
\usepackage{amsmath,amssymb}
\usepackage{verbatim}
\usepackage{dsfont,bm}
\usepackage{color}
\usepackage{graphicx}
\usepackage{amsthm}
\usepackage{enumerate}

\usepackage{xcolor}
\definecolor{darkgreen}{rgb}{0,0.75,0}
\definecolor{darkred}{rgb}{0.75,0,0}
\definecolor{darkmagenta}{rgb}{0.5,0,0.5}
\usepackage[colorlinks,citecolor=darkgreen,linkcolor=darkred,urlcolor=darkmagenta]{hyperref}
\newtheorem{theorem}{Theorem}[section]
\newtheorem{thm}[theorem]{Theorem}
\newtheorem{corollary}[theorem]{Corollary}
\newtheorem{cor}[theorem]{Corollary}
\newtheorem{lemma}[theorem]{Lemma}
\newtheorem{lem}[theorem]{Lemma}
\newtheorem{proposition}[theorem]{Proposition}
\newtheorem{prop}[theorem]{Proposition}
\newtheorem{definition}[theorem]{Definition}

\newtheorem{assumption}[theorem]{Assumption}
\newtheorem{remark}[theorem]{Remark}

\newtheorem{problem}[theorem]{Problem}
\numberwithin{equation}{section}

\def\be{\begin{equation}}
\def\ee{\end{equation}}
\def\bes{\begin{equation*}}
\def\ees{\end{equation*}}

\newcommand{\Cpc}[0]{\operatorname{Cap}}
\newcommand{\mr}[1]{{\tt \href{http://www.ams.org/mathscinet-getitem?mr=#1}{MR#1}}}
\newcommand{\arxiv}[1]{{\tt \href{http://arxiv.org/abs/#1}{arXiv:#1}}}
\newcommand{\set}[1]{\left\{ #1 \right\}}

\newcommand{\Sett}[2]{\left\{ #1  : \, #2 \right\}}

\newcommand{\BC}[0]{\hyperlink{bc}{$\operatorname{(MD)}$}}
\newcommand{\BG}[0]{\hyperlink{bg}{$\mathbf{\operatorname{(BG)}}$}}
\newcommand{\abs}[1]{{\left\vert\kern-0.25ex #1
    \kern-0.25ex\right\vert}}
\newcommand\norm[1]{\left\lVert#1\right\rVert} %norm
\newcommand{\one}{\mathds{1}} %indicator
\newcommand{\loc}[0]{\operatorname{loc}}
\newcommand{\diag}[0]{\operatorname{diag}}
\DeclareMathOperator*{\esssup}{ess\,sup}
\DeclareMathOperator*{\essinf}{ess\,inf}
\newcommand{\rst}{\operatorname{\mathbf{rst}}}
\newcommand{\ext}{\operatorname{\mathbf{ext}}}
\newcommand{\compl}{c}

\def\sA {{\mathcal A}}  \def\sC {{\mathcal C}}
\def\sD {{\mathcal D}} \def\sE {{\mathcal E}} \def\sF {{\mathcal F}}
  
  \def\sL {{\mathcal L}}
\def\sM {{\mathcal M}}

  \def\sX {{\mathcal X}}
\def\sY {{\mathcal Y}}

 \def\bE {{\mathbb E}} 
\def\bG {{\mathbb G}}  
  
 \def\bN {{\mathbb N}} 
\def\bP {{\mathbb P}}  \def\bR {{\mathbb R}}
  
\def\bV {{\mathbb V}}  
 \def\bZ {{\mathbb Z}}

\def\med{\medbreak\noindent}

\def\sms{\smallskip}
\def\ms{\medskip}

\def\sm{\smallskip\noindent}

\def\ignore#1{}

\def\ol{\overline}           

\def\lam {\lambda}  
\def\eps{\varepsilon}
\def\th{\theta} 
\def\vp{\varphi}
\def\Gam{\Gamma} \def\gam{\gamma}

\def\to {\rightarrow}

\def\pd {\partial}
\def\q{\quad} 
\def\dint{\int\kern-.6em\int}
\def\grad{\nabla}

\def\supp{\mathop{{\rm supp}}}

\def \half {{\tfrac12}}
\def\fract{\tfrac}

\def\wt{\widetilde}

\def\be{\begin{equation}}
\def\ee{\end{equation}}
\def\bes{\begin{equation*}}
\def\ees{\end{equation*}}
\def\ba{\begin{align}}
\def\ea{\end{align}}
\def\xxea{\end{align}}
\def\bas{\begin{align*}}
\def\eas{\end{align*}}

\def\nn{\nonumber} 

\def\proof{{\smallskip\noindent {\em Proof. }}}
\def\qed{{\hfill $\square$ \bigskip}}

\def\Ceff{C_{\rm eff}}

\begin{document}

\font\titlefont=cmbx14 scaled\magstep1
\title{\titlefont Stability of elliptic Harnack inequality} 

\author{
Martin T. Barlow\footnote{Research partially supported by NSERC (Canada)}, \,
Mathav Murugan\footnote{Research partially supported by NSERC (Canada) and the Pacific Institute for
 the Mathematical Sciences}
}

\maketitle

%\ver

\begin{abstract}
We prove that the elliptic Harnack inequality (on a manifold, graph, or suitably regular metric measure space)
 is stable under bounded perturbations, as well as rough isometries.


%\vskip.2cm \noindent {\it American Mathematical Society 
%Subject Classification:  Primary 35K08 Secondary 31B05 47D07}
\end{abstract}

\section{Introduction} \label{sec:intro}
A well known theorem of Moser \cite{Mo1} is that an elliptic Harnack 
inequality (EHI) holds for solutions associated with uniformly elliptic
divergence form PDE. Let $\sA$ be given by
\be
 \sA f(x) = \sum_{i,j=1}^d \frac{\pd }{\pd x_i} \Big( a_{ij}(x) \frac{\pd f}{\pd x_j} \Big),
\ee
where $(a_{ij}(x), x \in \bR^d)$ is bounded, measurable and uniformly elliptic. Let
$h$ be a non-negative $\sA$-harmonic function in a domain $B(x,2R)$, 
and let $B=B(x,R) \subset B(x,2R)$. Moser's theorem states that
there exists a constant $C_H$, depending only on $d$ and the ellipticity constant
of $a_{..}(\cdot)$, such that
\be \label{e:ehim1}
 \esssup_{B(x,R)} h \le C_H  \essinf_{B(x,R)} h.
\ee
A few years later Moser \cite{Mo2,Mo3} extended this to obtain
a parabolic Harnack inequality (PHI) for solutions $u=u(t,x)$ to the heat equation associated with $\sA$:
\be \label{e:Ahe}
 \frac{\pd u}{\pd t} = \sA u. 
\ee
This states that if $u$ is a non-negative solution to \eqref{e:Ahe} in
a space-time cylinder $Q=(0,T) \times B(x,2R)$, where $R=T^2$, then writing
$Q_- = (T/4,T/2) \times B(x,R)$, $Q_+= (3T/4,T) \times B(x,R)$,
\be \label{e:phim1}
 \esssup_{Q_-} u \le C_P  \essinf_{Q_+} u.
\ee
If $h$ is harmonic then $u(t,x) = h(x)$ is a solution to \eqref{e:Ahe}, so the PHI implies the EHI. 
The methods of Moser are very robust, and have been extended to manifolds,  
metric measure spaces, and graphs -- see \cite{BG,Sal92,St,De1,MS1}. 

The EHI and PHI have numerous applications, and in particular give
a priori regularity for  solutions to \eqref{e:Ahe}.
It is well known that Harnack inequality is useful beyond the linear elliptic and parabolic equations mentioned above. For instance variants of Harnack inequality apply to non-local operators, non-linear equations and geometric evolution equations including the Ricci flow and mean curvature flow -- see the survey \cite{Kas}.

 S.T. Yau and his collaborators \cite{Yau,CY,LY}
 developed a completely different approach to Harnack inequalities based on gradient estimates. 
 \cite{Yau} proves the Liouville property for Riemannian manifolds with non-negative Ricci curvature using gradient estimates for positive harmonic functions. 
 A local version of these gradient estimates was 
 given by Cheng and Yau in \cite{CY}. 
 Let $(M,g)$ a Riemannian manifold whose Ricci curvature is bounded below by $-K$ for some $K \ge 0$. Fix $\delta \in (0,1)$. 
Then there exists $C>0$, depending only on $\delta$ and $\dim(M)$, 
such that any positive solution $u$ of the Laplace equation 
$\Delta u=0$ in $B(x,2r) \subset M$ satisfies
 \[
 \abs{\nabla \ln(u)} \le C (r^{-1} + \sqrt{K}) \q \mbox{in $B(x,2 \delta r)$.}
 \]
 Integrating this estimate along geodesics immediately yields a local version of the EHI. 
 In particular, any $u$ above satisfies
 \[
 u(z)/u(y) \le \exp(C (1+\sqrt{K}r), \q z,y \in B(x,2\delta r).
 \]
For the case of manifolds with non-negative Ricci curvature we have $K=0$, 
and so obtain the EHI. This gradient estimate was extended to the parabolic
setting by Li and Yau \cite{LY}.
See \cite[p. 435]{Sal95} for a comparison between the gradient estimates of \cite{Yau,CY,LY} 
and the Harnack inequalties of Moser \cite{Mo1,Mo2}.

A major advance in understanding the PHI was made in 1992 by Grigoryan and Saloff-Coste
\cite{Gr0, Sal92}, who proved that the PHI is equivalent to two conditions: volume
 doubling (VD) and a family of Poincar\'e inequalities (PI).  The context
 of  \cite{Gr0, Sal92} is 
 the Laplace-Beltrami operator on Riemannian manifolds, but the basic
 equivalence VD+PI $\Leftrightarrow$ PHI also holds for graphs and metric measure
 spaces with a Dirichlet form -- see \cite{De1,St}.
This characterisation of the PHI implies that it is stable with respect to rough isometries -- see
\cite[Theorem 8.3]{CS}. 
For more details and a survey of the literature see the introduction of \cite{Sal95}.

One consequence of the EHI is the Liouville property -- that all bounded
harmonic functions are constant. However, the Liouville property is
not stable under rough isometries -- see \cite{Lyo}. 
See \cite[Section 5]{Sal04} for a survey of  related results and open questions.

Using the gradient estimate in \cite[Proposition 6]{CY}, 
Grigor'yan  \cite[p. 340]{Gr0} remarks that there exists a two dimensional Riemannian manifold that satisfies the EHI but does not satisfy the PHI.
In the late 1990s further examples inspired by analysis on fractals were given  -- see \cite{BB1}.
The essential idea behind the example in \cite{BB1} is that
if a space is roughly isometric to an infinite Sierpinski
carpet, then a PHI holds, but with anomalous space time scaling given by
$R= T^\beta\vee T^2$, where $\beta>2$. This PHI implies the EHI, but the
standard PHI (with $R=T^2$) cannot then hold. (One cannot have the
PHI with two asymptotically distinct space-time scaling relations.) 
\cite{BB3, BBK} prove that the anomalous PHI$(\Psi)$ with scaling 
$R=\Psi(T) = T^{\beta_1} \one_{(T \le 1)} + T^{\beta_2} \one_{(T >1)}$  is stable under rough isometries.
These papers also prove
that PHI$(\Psi)$ is equivalent to volume doubling, a family of Poincar\'e inequalities
with scaling $\Psi$, and a new inequality which controlled the energy
of cutoff functions in annuli, called a {\em cutoff Sobolev inequality}, and
denoted CS$(\Psi)$.
The papers \cite{BB3, BBK} proved the PHI by Moser's argument, but the more recent
papers \cite{AB, GHL}  use de Giorgi's argument and a mean value inequality to obtain
similar results, but with a
to simpler form of the cutoff Sobolev inequality. In addition, an important point 
for this paper, \cite{GHL} does not require the underlying metric space to be a length space.

A further example of weighted Laplace operators 
on  Riemannian manifolds that satisfy EHI but not PHI is given in \cite[Example 6.14]{GS}.
Consider the second order differential operators $L_\alpha$ on $\bR^n, n \ge 2$ given by
\[
L_\alpha= \left(1 + \abs{x}^2\right)^{-\alpha/2} 
\sum_{i=1}^n \frac{\partial}{\partial x_i} \left(\left(1 + \abs{x}^2\right)^{\alpha/2} \frac{\partial}{\partial x_i} \right)
= \Delta  +  \alpha\frac{ x. \nabla}{1 + \abs{x}^2}.
\]
Then $L_\alpha$ satisfies the PHI if and only if $\alpha> -n$ but satisfies 
the EHI for all $\alpha \in \bR$. Weighted Laplace operators of this kind
arise naturally in the context of Schr\"odinger operators and conformal transformations of Riemannian metrics -- see \cite[Section 6.4 and 10]{Gri06}. 

These papers left open the problem of the stability of the EHI, and also the question 
of finding a satisfactory characterisation of the EHI. This problem is mentioned
in  \cite{Gri95}, \cite[Question 12]{Sal04} and \cite{Kum}. 
In \cite{GHL0}, the authors write  ``An interesting (and obviously hard) question is the characterization of the elliptic
Harnack inequality  in more geometric terms -- so far nothing is known, not even a conjecture."

In \cite{De2}, Delmotte gave an example of a graph which satisfies the EHI but for which
(VD) fails; his example was to take the join of the infinite Sierpinski gasket
graph with another (suitably chosen) graph. 
 This example shows that any attempt to characterize the EHI
must tackle the difficulty that different parts of the space may have 
different space-time scaling functions. Considerable progress
on this was made by R. Bass \cite{Bas}, but his result requires 
volume doubling, as well as some additional hypotheses on capacity. 

As Bass remarks, all the robust proofs of the EHI, using the methods of 
De Giorgi, Nash or Moser, use the 
volume doubling property in an essential way, as well as 
Sobolev and Poincar\'e type inequalities. The starting point for this paper is the observation that 
a change of the symmetric measure (or equivalently a time change of the process) does not affect the 
sheaf of harmonic functions  on bounded open sets. On the other hand 
properties such as volume doubling or Poincar\'e inequality are
not in general preserved by this transformation.

Conversely,  given a space satisfying the EHI, one could seek to construct a `good' measure $\mu$ such
that volume doubling, as well as additional Poincar\'e and Sobolev inequalities
do hold with respect to $\mu$; this is indeed the approach of this paper.
Our main result, Theorem \ref{T:eeprime}, is that the EHI is stable. Our methods also give 
a characterization of the EHI by properties that are easily seen to be stable under 
perturbations -- see Theorem \ref{T:main-new}.

Our main interest is the EHI for manifolds and graphs. 
To handle both cases at once we 
work in the general context of metric measure spaces. 
So we consider a complete, locally compact, separable, geodesic (or length) 
metric space $(\sX,d)$ with a Radon measure $m$ which has full support,
so that $m(U)>0$ for all non-empty, open $U$. We call this a {\em metric measure space}.
Let $(\sE,\sF^m)$ be a strongly local Dirichlet form on $L^2(\sX,m)$ -- see \cite{FOT}.
We call the quintuple $(\sX,d,m,\sE,\sF^m)$ a {\em measure metric space
with Dirichlet form}, or {\em MMD space}.
We write $B(x,r) =\{y: d(x,y)<r\}$ for open balls in $\sX$, and 
given a ball $B=B(x,r)$ we sometimes use the notation
$\theta B$ to denote the ball $B(x,\theta r)$.
We assume  $(\sX,d)$ has infinite radius, so that 
$\sX - B(x,R) \neq \emptyset$  for all $R>0$. 
See Section \ref{sec:mmd} for more details of these spaces, and the definitions
of harmonic functions and capacities in this context.

Our two fundamental examples are Riemannian manifolds and the cable systems of graphs.
If $(\sM,g)$ is a Riemannian manifold we take $d$ and $m$ to be
the Riemannian distance and measure respectively, and define the Dirichlet 
form to be the closure of the symmetric bilinear form
\bes
 \sE(f,f) = \int_\sX |\grad_g f|^2 dm, \q f \in C^\infty_0(\sM).
\ees

Given a graph $\bG=(\bV,E)$ the {\em cable system} of $\bG$ is the metric space obtained
by replacing each edge by a copy of the unit interval, glued together in the obvious way.
For a graph with uniformly bounded vertex degree the EHI 
for the graph is equivalent to the EHI for its cable system, and so
our theorem also implies stability of the EHI for graphs.
See Section \ref{sec:ex} for more details of both these examples.

Since our main spaces of interest are regular at small length scales, we will  avoid a number of
technical issues which could arise for general MMD spaces by making 
two assumptions of local regularity:
Assumptions \ref{A:green} and \ref{A:BG}. Both our
main examples satisfy these assumptions -- see Section \ref{sec:ex}.

The hypothesis of volume doubling plays an important role
in the study of heat kernel bounds for the process $X$, and as mentioned
above is a necessary condition for the PHI.

\begin{definition}[Volume doubling property]
{\rm We say that a Borel measure $\mu$ on a metric space $(\sX,d)$ 
satisfies the \emph{volume doubling property}, if $\mu$ is non-zero and there 
exists a constant $C_V <\infty$ such that
\be
\mu(B(x,2r)) \le  C_V \mu(B(x,r))
\ee
for all $x \in \sX$ and for all $r>0$.
}\end{definition}

\begin{definition}\label{D:ehi}
{\rm We say that $(\sX,d,m,\sE,\sF^m)$ satisfies
the {\em elliptic Harnack inequality} (EHI) if there exist
constants $1<A,C_H<\infty$ such that for any $x \in \sX$ and $R>0$, 
for any nonnegative harmonic function $h$  on a ball $B(x,AR)$ one has
\be \label{e:ehi2}
   \esssup_{B(x,R)} h \le C_H  \essinf_{B(x,R)} h.
\ee
}\end{definition}
If $(\sX,d)$ is a geodesic metric space and the above inequality holds for 
some value of $A>1$, then it holds for any other $A'>1$ with a constant $C_H(A')$. 
If the EHI holds, then iterating the condition \eqref{e:ehi2} gives a.e.~H\"older continuity
of harmonic functions, and it follows that any harmonic function has a continuous modification.

Our first main theorem is
\begin{thm} \label{T:eeprime}
Let $(\sX,d,m)$ be a length metric measure space, and $(\sE,\sF)$
be a strongly local Dirichlet form on $L^2(\sX,m)$. Suppose that
Assumptions \ref{A:green} and \ref{A:BG} hold.  Let
$(\sE',\sF)$ be a strongly local Dirichlet form on $L^2(\sX,m')$
which is equivalent to $\sE$, so that there exists $C<\infty$ such that
\begin{align*}
 C^{-1} \sE(f,f) &\le \sE'(f,f) \le  C \sE(f,f) \q \hbox{ for all } f \in \sF,\\
  C^{-1} m(A) &\le m'(A) \le  C m(A) \q \hbox{ for all measurable sets } A.
\end{align*}
Suppose that $(\sX,d,m,\sE,\sF)$ satisfies the elliptic Harnack inequality.
Then the  EHI holds for $(\sX,d,m',\sE',\sF)$. 
\end{thm}

We now state some
consequences of Theorem \ref{T:eeprime} for Riemannian manifolds and graphs.
We say that two Riemannian manifolds $(M,g)$ and $(M',g')$ are 
{\em quasi isometric} if there exists a diffeomorphism $\phi:(M,g) \to (M',g')$ 
and a constant $K \ge 1$ such that 
\[
K^{-1} g(\xi,\xi) \le g'( d\phi(\xi),d\phi(\xi)) \le 	K g(\xi,\xi), \q \mbox{ for all } \xi \in TM. 
\]

Let $(M,g)$ be a Riemannian manifold and let $\operatorname{Sym}(TM)$ denote 
the bundle of symmetric endomorphisms of the tangent bundle $TM$. 
We say that $\sA$ is an \emph{uniformly elliptic operator in divergence form}
if there exists  $A:M \to \operatorname{Sym}(TM)$ a measurable section of 
$\operatorname{Sym}(TM)$ and a constant $K \ge 1$ such that 
\bes 
	K^{-1} g(\xi,\xi) \le g(A\xi,\xi) \le K g(\xi,\xi), \q \forall \xi \in TM,
\ees   
such that $\sA (\cdot)= \operatorname{div} \left( A \nabla (\cdot) \right)$. 
Here $\operatorname{div}$ and $\nabla$ denote the Riemannian divergence 
and gradient respectively.

\begin{theorem} \label{T:mainconsq}
(a) Let $(M,g)$ be a Riemannian manifold that is quasi isometric
to a manifold whose Ricci curvature is bounded below, 
and let $\Delta$ denote the corresponding Laplace-Beltrami operator. 
If $(M,g)$ satisfies the EHI for non-negative solutions of $\Delta u=0$, 
then it satisfies the EHI for  non-negative  solutions of $\sA u=0$, 
where $\sA$ is any uniformly elliptic operator in divergence form. \\
(b)  Let $(M,g)$ and $(M',g')$ be two Riemannian manifolds 
that are quasi isometric to a manifold whose Ricci curvature is bounded below. 
Let $\Delta$ and $\Delta'$ denote the corresponding Laplace-Beltrami operators.
Then,  non-negative  $\Delta$-harmonic functions satisfy the 
EHI, if and only if non-negative  ${\Delta'}$-harmonic functions satisfy the EHI.
\end{theorem}

\begin{theorem} \label{T:mainconsqG}
Let $\bG=(\bV,E)$ and $\bG'=(\bV',E')$ be bounded degree graphs, which are roughly isometric.
Then the EHI holds for $\bG'$ if and only if it holds for $\bG$.
\end{theorem}

\begin{remark}
{\rm 
(1) Theorem \ref{T:mainconsq}(a) is a generalization 
of Moser's elliptic Harnack inequality \cite{Mo1}.
The parabolic versions of  (a) and (b) are due to \cite{Sal92b}.
For (b) note the the manifold $(M,g)$ might not have Ricci curvature bounded 
below and hence the methods of \cite{Yau,CY} will not apply. 
A parabolic version of Theorem \ref{T:mainconsqG} is essentially due to \cite{De1}.\\
(2) As proved in \cite{Lyo}, the Liouville property is not stable under rough isometries.
}\end{remark}

The outline of our argument is as follows. In Section \ref{sec:conEHI}
using the tools of potential theory 
we prove that the EHI implies certain regularity properties 
for Green's functions and capacities. 
The main result of this section (Theorem \ref{T:BC}) is that the EHI implies that
$(\sX,d)$ has the metric doubling property.

\begin{definition} \label{e:MDdef2}
{\rm 
The space $(\sX,d)$ satisfies the {\em metric doubling property} 
\hypertarget{bc}{$\operatorname{(MD)}$} if there exists $M<\infty$ such that
any ball $B(x,R)$ can be covered by $M$ balls of radius $R/2$. 
}\end{definition}

An equivalent definition is that there exists $M'<\infty$ such that 
any ball $B(x,R)$ contains at most $M'$ points which are all a distance
of at least $R/2$ from each other. %R/2?
We will frequently use the fact that (MD) 
holds for $(\sX,d)$ if and only if $(\sX,d)$ has finite Assouad dimension. 
Recall that  the Assouad dimension 
is the infimum of all numbers $\beta >0$ with the property that every ball of radius $r>0$ has at most $C \eps^{-\beta}$ disjoint points of mutual distance at  least $\eps r$ for some $C \ge 1$ independent of the ball. (See \cite[Exercise 10.17]{Hei}.) Equivalently, this is  the infimum of all numbers $\beta >0$ with the property that 
every ball of radius $r>0$ can be covered by at most $C \eps^{-\beta}$ balls of 
radius $\eps r$ for some $C \ge 1$ independent of the ball. 

It is well known that volume doubling implies metric doubling.
A partial converse also holds: 
if  $(\sX,d)$ satisfies \BC\hspace{0.3mm} then there exists a Radon
measure $\mu$ on $\sX$ such that  $(\sX,d,\mu)$ satisfies (VD).
This is a classical result
due to Vol'berg and Konyagin \cite{VK} in the case of compact spaces, and
Luukkainen and Saksman \cite{LuS} in the case of general complete spaces. 
For other proofs see  \cite{Wu} and \cite[Chapter 13]{Hei}), and also
 \cite[Chapter 10]{Hei} for a survey of some conditions equivalent to (MD).

The measures constructed in these papers are very far from being unique. 
In Section \ref{sec:goodVD},
using the approach of  \cite{VK}, we show that if $\sX$ satisfies the EHI and
Assumptions \ref{A:green} and \ref{A:BG} then we can
construct a `good' doubling measure $\mu$ which is absolutely continuous with respect
to $m$,  and connects capacities with the measures of balls  in a suitable
fashion -- see Definition \ref{d:goodmeas} and Theorem \ref{T:meas}.

At this point we could use some extensions of the methods of \cite{Bas, GHL} to prove the stability of the
EHI. However, a quicker approach, which we follow in Section \ref{sec:qs}, is to use ideas 
from the theory of quasisymmetric transformations of metric spaces. (See \cite{Hei} for
an introduction to this theory, and \cite{Ki2} for applications 
to heat kernels.) 
These transformations do not distort annuli too much, and therefore preserve the EHI -- see Lemma \ref{l:ehiqs}. 
In Section \ref{sec:qs} we prove that there exists a new metric $d_\Psi$ on $\sX$ and a constant
$\beta>0$ such that the new space $(\sX, d_\psi, \mu)$ satisfies 
Poincar\'e and cutoff energy inequalities with respect to a global space-time scaling 
relation of the form $R = T^\beta$ -- see Theorem \ref{T:PI-CS}. These inequalities are
stable with respect to bounded perturbations of the Dirichlet form.  
While metric $d_\Psi$ is not geodesic the main theorem of 
\cite{GHL} does apply in this context, and gives that 
the PI and CS inequalities  in Theorem \ref{T:PI-CS} imply the EHI. This 
gives the stability of the EHI, as well as a stable characterization -- see Theorem 
\ref{T:main-new}.

In Section \ref{sec:ex} we return to our two main classes of examples, weighted Riemannian
manifolds and weighted graphs. We show that they both satisfy our local regularity
hypotheses Assumptions  \ref{A:green} and \ref{A:BG}, and give the (short)
proof of Theorem \ref{T:mainconsq}.

The final Section \ref{sec:qi} formulates the class of rough isometries which
we consider, and states our result on the stability of the EHI under rough isometries.
Since rough isometries only relate spaces at large scales,
and the EHI is a statement which holds at all length scales,  
any statement of stability under rough isometries requires that the family
of spaces under consideration satisfies suitable local regularity hypotheses.

A characterization of the EHI in terms of effective resistance
(equivalently capacity) was suggested in \cite{B1}. G. Kozma \cite{Ko}
gave an illuminating counterexample -- a spherically symmetric tree. This example does
not satisfy (MD), and at the end of Section 7 we suggest a modified 
characterization, which is the `dumbbell condition' of \cite{B1} together with (MD).

We use $c, c', C, C'$ for strictly positive constants, which may change value 
from line to line. Constants with numerical subscripts will keep the same value in each
argument, while those with letter subscripts will be regarded as constant
throughout the paper.
The notation $C_0 = C_0(a,b)$ means that the constant
$C_0$ depends only on the constants $a$ and $b$.

\section{Metric measure spaces with Dirichlet form} \label{sec:mmd}

In this section give some background on MMD spaces, and give our two
assumptions of local regularity.
We take $(\sX,d)$ to be a locally compact metric space
with infinite radius, and $m$ to be a Radon measure on $(\sX,d)$ with full support.
Let $(\sE,\sF^m)$ be a strongly local Dirichlet form on $L^2(\sX,m)$.
We call  $(\sX,d,m,\sE,\sF^m)$ a {\em measure metric space
with Dirichlet form}, or {\em MMD space}.
Except in Section \ref{sec:qs} we will assume also that $(\sX,d)$ is a length space.

In the context of MMD spaces Poincar\'e and Sobolev inequalities 
involve integrals with respect to the  {\em energy measures} $d\Gam(f,f)$ -- formally
these can be regarded as $|\grad f|^2 dm$.
For bounded  $f \in \sF^m$ the measure $d\Gam(f,f)$
is  defined to be the unique measure 
such that for all bounded $g \in \sF^m$ we have
$$ \int g d\Gam(f,f) = 2 \sE(f,fg) - \sE(f^2,g). $$
We have
$$ \sE(f,f) = \int_{\sX} d\Gam(f,f). $$
For a Riemannian manifold
$d\Gam(f,f) =  |\grad_g f|^2 dm. $

Associated with $(\sE, \sF^m)$ is a semigroup $(P_t)$ and its infinitesimal
generator $(\sL, \sD(\sL))$. The operator $\sL$ satisfies
\be
  - \int( f \sL  g) dm = \sE(f,g), \q f \in \sF^m, g \in \sD(\sL);
\ee
in the case of a Riemannian manifold $\sL$ is the Laplace-Beltrami operator.
$(P_t)$ is the semigroup of a continuous Hunt process
$X=(X_t, t \in [0,\infty), \bP^x, x \in \sX)$.

We define capacities for $(\sX,d,m,\sE,\sF^m)$ as follows. For 
a non-empty open subset $D \subset \sX$, let $\sC_0(D)$ denote the space of all continuous functions with compact support in $D$. Let $\sF_D$ denote the closure of $\sF^m \cap \sC_0(D)$ with 
respect to the $\sqrt{ \sE(\cdot,\cdot)+ \norm{\cdot}_2^2}$-norm. 
By $A \Subset D$, we mean that the closure of $A$ is a compact subset of $D$.
For $A \Subset D$ we set
\be \label{e:capdef}
 \Cpc_D(A) = \inf\{ \sE(f,f): f \in \sF_D \mbox{ and  $f \ge 1$ in a neighbourhood of $A$} \}.
\ee

It is clear from the definition that if $A_1 \subset A_2 
\Subset D_1 \subset D_2$ then
\be \label{e:capmon}
 \Cpc_{D_2}(A_1) \le \Cpc_{D_1}(A_2). 
\ee
We can consider $\Cpc_D(A)$ to be the effective conductance 
between the sets $A$ and $D^\compl$ if we regard $\sX$ as an electrical
network and $\sE(f,f)$ as the energy of the function $f$. 
A statement depending on $x \in B$ is said to hold quasi-everywhere on $B$ (abbreviated as q.e.~on $B$), if there exists a set $N \subset B$ of zero capacity such that the statement if true for every $x \in B \setminus N$.
It is known that $(\sE,\sF_D)$ is a regular Dirichlet form on $L^2(D,m)$ and 
\be
 \sF_D = \{ f \in \sF^m: \wt{f}=0\, \mbox{ q.e.~on $D^\compl$}\},
\ee
where $\wt{f}$ is any quasi continuous representative of $f$ (see \cite[Corollary 2.3.1 and Theorem 4.4.3]{FOT}).
Functions in the extended Dirichlet space will always represented by 
their quasi continuous version (cf. \cite[Theorem 2.1.7]{FOT}), 
so that expressions like $\int f^2 d\Gamma(\vp,\vp)$ 
are well defined.

Given an open set $U \subset \sX$, we set 
\begin{align*}
\sF_{\loc}(U) = \{ h \in &L^2_{\loc}(U): \mbox{for all relatively compact } V \subset U, \\
& \mbox{ there exists } h^\# \in \sF^m, \mbox{ s.t. } h \one_V = h^\# \one_V  \, m \mbox{-a.e.}\}.
\end{align*}

\begin{definition}
{\rm A function $h:U \to \bR$ is said to be \emph{harmonic} in an open set $U \subset \sX$, if $h \in \sF_{\loc}(U)$ and satisfies $\sE(f,h)=0$ for all $f \in \sF^m \cap \sC_0(U)$. Here $\sE(f,h)$ can be unambiguously defined as $\sE(f,h^\#)$ where  $h = h^\#$ in a precompact open set containing $\operatorname{supp}(f)$ and  $h^\# \in \sF^m$.
}\end{definition}

This definition implies that $\sL h=0$ in $D$ provided that $h$ is 
in the domain of $\sL_D$.

Next, we define the Green's operator and Green's function. 

\begin{definition}
{\rm Let $D$ be a bounded open subset of $\sX$.
Let $\sL_D$ denote the generator of the Dirichlet form $(\sE,\sF_D,L^2(D,m))$ 
and assume that
\be \label{e:lammindef}
\lambda_{\min}(D)= \inf_{f \in \sF_D \setminus \set{0}} \frac{\sE(f,f)}{\norm{f}_2^2} > 0.
\ee
We define the inverse of $-\sL_D$ as the \emph{Green operator} $G_D=(-\sL_D)^{-1}: L^2(D,m) \to L^2(D,m)$. We say a jointly measurable function 
$g_D(\cdot,\cdot): D \times D \to \bR$ is the \emph{Green function} for $D$ if 
\[
G_Df(x)= \int_D g_D(x,y) f(y) \, m(dy) \,\, \mbox{ for all } f \in L^2(D,m) \mbox{ and for $m$ a.e. $x \in D$.} 
\]
} \end{definition}

\begin{assumption}  \label{A:green}
{\rm (Existence of Green function) 
For any bounded, non-empty open set $D \subset \sX$, we 
assume that 
$\lambda_{\min}(D)>0$ and that there exists a Green function $g_D(x,y)$  for $D$
defined for $(x,y) \in D\times D$ with the following properties:
 \begin{enumerate}[(i)]
 \item (Symmetry) $g_D(x,y)= g_D(y,x) \ge 0$ for all $(x,y) \in D \times D \setminus \diag$;
 \item (Continuity) $g_D(x,y)$ is jointly continuous in $(x,y) \in D \times D \setminus \diag$;
 \item (Maximum principles) If $x_0 \in U \Subset D$, then
 \begin{align*}
 \inf_{U \setminus \set{x_0}} g_D(x_0,\cdot) = \inf_{\partial U} g_D(x_0,\cdot),
\qquad \sup_{D \setminus U} g_D(x_0, \cdot) = \sup_{\partial U}g_D(x_0,\cdot).
 \end{align*}
 \item (Harmonic) For any fixed $x \in D$, the function $y \mapsto g_D(x,y)$ is 
 in $\sF_{\loc}(D \setminus \set{x})$ and is harmonic in $D \setminus \set{x}$.
 \end{enumerate}
 Here $\diag$ denotes the diagonal in $D \times D$.
}\end{assumption}

\begin{remark} \label{R:ellipm}
{\rm 
Note that changing the measure $m$ to an 
equivalent Radon measure $m'$ does not affect either the the capacity of bounded sets
or the class of harmonic functions. Further, if $f_1, f_2 \in C(\sX) \cap \sF_D$  then
writing $\langle \cdot , \cdot \rangle_m$ for the inner product in $L^2(m)$,
\be
 \sE( G_D f_1, f_2) = \langle f_1, f_2 \rangle_m,
\ee
and it follows that $g_D(\cdot,\cdot)$ is also not affected by this change of measure.
} \end{remark}

Our second key local regularity assumption is as follows.

\begin{assumption} \label{A:BG}
{\rm (Bounded geometry or \hypertarget{bg}{($\mathbf{\operatorname{BG}}$)}).
We say that a MMD space $(\sX,d,m,\sE,\sF^m)$ satisfies (BG) if 
there exist $r_0 \in (0, \infty]$ and $C_L < \infty$ such that the following hold:
\begin{enumerate}
\item[(i)] (Volume doubling property  at small scales). 
For all $x \in \sX$ and for all $r \in (0,r_0]$ we have
	\begin{equation} \label{e:vds}
	\frac{m(B(x,2r))}{m(B(x,r))} \le C_L.
	\end{equation}
\item[(ii)] (Expected occupation time growth at small scales.) 
There exists $\gamma_2>0$ such that for 
all $x_0 \in \sX$ and for all $0<s \le r  \le r_0$ we have
	\begin{equation}\label{e:eois}
	\frac{ m(B(x,s))} { \Cpc_{B(x,8s)}(B(x,s))}
	\frac{ \Cpc_{B(x,8r)}(B(x,r))} {m(B(x,r))} 
	\le C_L \left( \frac{s}{r}\right)^{\gamma_2} .
	\end{equation}
\end{enumerate}
}\end{assumption}

See Section \ref{sec:ex} for the verification of 
Assumptions \ref{A:BG} and \ref{A:green} for our two
main cases of interest, weighted Riemannian manifolds with
Ricci curvature bounded below, and the cable system
of graphs with uniformly bounded vertex degree.
The condition (BG)  is a robust one, because under mild conditions it is
preserved under bounded perturbation of conductance in a weighted graph,
 and quasi isometries of weighted manifolds.

\section{Consequences of EHI} \label{sec:conEHI}

\sm
Throughout this section we assume that $(\sX,d,m,\sE,\sF^m)$ 
is a MMD space which 
satisfies Assumption \ref{A:green}, as well as the EHI with constant $C_H$.
In addition we assume that $(\sX,d)$ is a length space,
and will  write $\gam(x,y)$ for a geodesic between $x$ and $y$.
Recall that $(X_t)$ is the Hunt process associated 
with $(\sE,\sF^m)$, and write for $F \subset \sX$,
\be \label{e:reg}
 T_F = \inf\{ t \ge 0: X_t \in F\}, \, \tau_F = T_{F^\compl}. 
 \ee

\begin{thm} \label{T:behit2}
Let $(\sX,d)$ satisfy the EHI. Then there exists a constant
$C_G = C_G(C_H)$ such that if $B(x_0,2R) \subset D$ then
\be
  g_D(x_0,y) \le C_G g_D(x_0,z) \hbox{ if  } d(x_0,y)=d(x_0,z)=R.
\ee
\end{thm}

\proof The proof of \cite[Theorem 2]{B1} carries over to this 
situation with essentially no change.
(In fact it is slightly simpler, since there is no need to make corrections
at small length scales.)
Note that since $g_D(\cdot, \cdot)$ is continuous off the diagonal, we can use the EHI
with $\sup$ and $\inf$ instead of $\esssup$ and $\essinf$.
\qed

\begin{cor} \label{C:hg3}
Let $B(x_0,2R) \subset D$.
Let $A \ge 2$. Then there exists a constant $C_1=C_1(C_H,A)$ such that
$$ g_D(x_0,x) \le C_1 g_D(x_0,y), \q \hbox{ for } x,y \in B(x_0,R)\setminus B(x_0, R/A). $$
\end{cor}

\proof We can assume $d(x,x_0) \ge d(x_0,y)$. Let $z$ be the point 
on $\gam(x_0,x)$ with $d(x_0,z)=d(x_0,y)$. Then we can compare
$g_D(x_0,y)$ and $g_D(x_0,z)$ by Theorem \ref{T:behit2}, and $g_D(x_0,z)$
and $g_D(x_0,x)$ by using a chain of balls with centres in $\gam(z,x)$. 
(The number of balls needed will depend on $A$.)
\qed

\begin{lemma} \label{L:greenc1}
Let $x_0 \in \sX$, $R>0$ and let $B(x_0,2R)\subset D$. 
There exists a constant
$C_0=C_0(C_H)$ such that if $x_1,x_2, y_1, y_2 \in B(x_0, R)$ with
$d(x_j,y_j) \ge R/4$ then
\be
   g_D(x_1, y_1) \le C_0 g_D(x_2,y_2).
\ee 
\end{lemma}

\proof 
A counting argument shows there exists a ball $B(z,R/9)\subset B(x_0,R)$ 
which contains none of the points $x_1,x_2, y_1, y_2$. Using
Corollary \ref{C:hg3} we have $g_D(x_1,y_1) \le c g_D(z,x_1)$, 
$g_D(z,x_1) \le c g_D(z,x_2)$, and $g_D(z,x_2) \le g_D(x_2,y_2)$,
and combining these comparisons gives the required bound. \qed

\begin{definition}
Set
\be
  g_D(x,r) = \inf_{y: d(x,y)=r} g_D(x,y). 
\ee
\end{definition}

The maximum principle implies that $g_D(x,r)$ is non-increasing in $r$.
An easy argument gives that if
$d(x,y) = r$ and $B(x,2r)\cup B(y,2r) \subset D$ then
\be
 g_D(x,r) \le C_G g_D(y,r). 
\ee

Let $D$ be a bounded domain in $\sX$, 
$A$ be Borel set, $A \Subset D \subset \sX$, and 
recall from \eqref{e:capdef} the definition of $\Cpc_D(A)$.
By  \cite[Theorem 4.3.3]{FOT}, \cite[Proposition A.2]{GH} 
there exists a function $h_{A,D} \in \sF_D$ called the \emph{equilibrium potential} 
such that $h_{A,D}=1$ q.e. in $A$ and $\sE(h_{A,D},h_{A,D})= \Cpc_{D}(A)$. 
The function $h_{A,D}(\cdot)$ is the hitting probability of the set $A$:
\be \label{e:cap-hit}
h_{A,D}(x)= \bP^x( T_A < \tau_D), \q \mbox{ for $x \in D$ quasi everywhere.} 
\ee
Further 
\be \label{e:cap-cutoff}
h_{A,D}(x)=1, \q \mbox{ quasi everywhere on $A$.} 
\ee
There exists a Radon measure $\nu_{A,D}$ called the {\em capacitary measure} or {\em equilibrium measure} that does not charge any set of zero capacity, supported on $\pd A$ such that 
$\nu_{A,D}(\partial A) = \Cpc_D(A)$ and  satisfies (cf. \cite[Lemma 2.2.10 and Theorem 2.2.5]{FOT} and \cite[Lemma 6.5]{GH})
\begin{align} \label{e:capmeas1}
\sE(h_{A,D}, v) &= \int_{\partial A} \wt{v} \, d\nu_{A,D}= \int_{D} \wt{v} \, d\nu_{A,D}, \q \mbox{for all $v \in \sF_D$;}\\\label{e:capmeas2}
 h_{A,D}(y) &= \nu_{A,D} G_D(y) = \int_{\pd A} \nu_{A,D}(dx) g_D(x,y), \q \mbox{for all $y \in D \setminus \pd A$.} 
\end{align}
Here $\wt{v}$ in \eqref{e:capmeas1} denotes a quasi continuous version of $v$.
By \cite[Theorem 2.1.5 and p.71]{FOT}
$\Cpc_D(A)$ can be expressed as 
\be \label{e:cap-char}
 \Cpc_D(A) = \inf\{ \sE(f,f): f \in \sF_D, f \ge 1 \mbox{ quasi everywhere on } A \}.
 \ee

\begin{lemma} \label{L:cap1}
Let $B(x_0,2r) \subset D$. Then
\be
 g_D(x_0,r) \le     \textstyle \Cpc_D(B(x_0,r))^{-1}  \le  C_G  g_D(x_0,r).
\ee
\end{lemma}

\proof Let $\nu$ be the capacitary measure for $B(x_0,r)$ with respect to $G_D$. Then
$\nu$ is supported by $\pd B(x,r)$ and by \eqref{e:capmeas2}
$$ 1 = \nu G(x_0) = \int_{\pd B}  g_D(x_0,z) \nu(dz). $$
Hence 
\bes
 \nu(B(x_0,r))  g_D(x_0,r) \le 1 \le \nu(B(x_0,r)) \sup_{z \in \pd B} g_D(x_0,z) \le C_G \nu(B(x_0,r)) g_D(x_0,r). 
 \ees
\qed

\begin{remark} \label{R:A2}
{\rm 
The assumption $B(x_0,2r) \subset D$ in 
Theorem \ref{T:behit2}, Corollary \ref{C:hg3}, Lemmas \ref{L:greenc1} and 
\ref{L:cap1} can be replaced with the assumption $B(x_0,Kr) \subset D$ for any fixed $K>1$.
}\end{remark}

\begin{lemma} \label{L:hit}
Let $B=B(x_0,R) \subset \sX$, and let $x_1  \in B(x_0, R/2)$,
$B_1= B(x_1, R/4)$. There exists $p_0 =p_0(C_H)$ such that 
\be \label{e:h11}
 \bP^y( T_{B_1} < \tau_B ) \ge p_0> 0 \hbox{ for } y \in B(x_0, 7R/8).
\ee
\end{lemma}

\proof
Let  $\nu$ be the capacitary measure for $B_1$ with respect to $G_B$,
and $h(x) = \nu G_B(x)=  \bP^x( T_{B_1} < \tau_B )$. 
Then $h$ is 1 on $B_1$, so by the maximum
principle it is enough to prove \eqref{e:h11} for $y \in B(x_0, 7R/8)$ with
$d(y, B_1) \ge R/16$. 

By Corollary \ref{C:hg3} (applied in a chain of balls if necessary) there exists
$p_0>0$ depending only on $C_H$ such that 
$g_B(y,z) \ge p_0 g_B(x_1,z)$ for $z \in \pd B_1$. 
Thus
$$ h(y)  \ge p_0   \int_{\pd B_1} g_B(x_1,z) \nu(dz) 
= p_0 \nu G_B(x_1) = p_0. $$
\qed

\begin{corollary}
Let $B(x_0,2R) \subset D$.
Then there exists $\th=\th(C_H)>0$ such that
if $0<s<r< R/2$ and $x \in B(x_0,R)$ then 
\be \label{e:ggrowth}
 \frac{ g_D(x,r) }{ g_D(x,s) } \ge c \Big( \frac{s}{r} \Big)^\th. 
 \ee
\end{corollary}

\proof Let $w \in \pd B(x,2s)$ and let $z \in \gam(x,w) \cap \pd B(x,s)$.
Applying the EHI on a chain of balls on $\gam(z,w)$ gives
$g_D(x,w) \ge c_1 g_D(x,z)$, and it follows that
$$ g_D(x,s) \le C_1 g_D(x,2s). $$
Iterating this estimate then gives \eqref{e:ggrowth} with $\th=\log_2 C_1$.
\qed

\begin{remark} \label{R:r2}
{\rm The example of $\bR^2$ shows that we cannot expect a corresponding
upper bound on $g_D(x,r)/g_D(x,s)$.
}\end{remark} 
 
 The key estimate in this Section is the following geometric consequence of the
 EHI. A weaker result proved with some of the same ideas, and in the graph
 case only, is given in \cite[Theorem 1]{B1}.
 
\begin{lemma} \label{L:l2}
Let $B=B(x_0,R) \subset \sX$. Let $\lam \in [\fract14, 1]$, 
$0< \delta \le 1/32$ and let $B_i = B(z_i,  \delta R)$, $i=1, \dots, n$ satisfy: \\
(1) $B_i \cap \pd B(x_0, \lam R) \neq \emptyset$, \\
(2) $B_i^* = B(z_i, 8 \delta R)$ are disjoint.\\
Then there exists a constant $C_1=C_1(C_H, \delta)$  such that $n \le C_1$.
\end{lemma}
 
\proof 
Let $y_i$ and $w_i$ be points on $\gam(x_0,z_i)$ with
$d(z_i,y_i)=3\delta R$ and $d(z_i,w_i)=5 \delta R$.
Let $A_i = \ol B(y_i, \delta)$.  By Lemma \ref{L:hit}
\be \label{e:h1}
\bP^x(  T_{ B_i} < \tau_{B_i^*} ) \ge p_1 >0, \; \mbox{ for all } x \in A_i. 
\ee

Now let $D= B - \cup_i B_i$, and let $N$ be the number of distinct balls
$A_i$ hit by $X$ before $\tau_D$. 
Write $S_1 < S_2 < \dots < S_N$ for the hitting times of these balls.
Let $G_{k i} =\{ N \ge k, X_{S_k} \in A_i \}$. 
On the event $G_{k i}$ if $X$ then hits $B_i$ before leaving $B_i^*$
we will have $N=k$. So using \eqref{e:h1}, for $k \ge 1$,
$$ \bP^x(N=k | N \ge k ) \ge p_1, $$
and thus $N$ is dominated by a geometric random variable with mean $1/p_1$.
Hence,
\be \label{e:h2}
\bE^{x_0} N \le 1/p_1. 
\ee
Now set
$$ h_i (x) = \bP^x( T_{A_i} < \tau_D ). $$
Then $h_i(y_i)=1$ and by Lemma \ref{L:hit} $h_i(w_i) \ge p_1$.
Using the EHI in a chain of balls we have
$h_i(x_0) \ge p_2=p_2(\delta)>0$.
Thus
$$ p_1^{-1} \ge \bE^{x_0} N = \sum_{i=1}^n h_i(x_0)
 \ge n p_2, $$
 which gives an upper bound for $n$. \qed

\begin{thm} \label{T:BC}
Let $(\sX,d,m,\sE,\sF^m)$ satisfy EHI. Then $(\sX,d)$ satisfies the metric doubling property \BC.
\end{thm}

\proof Let $\delta = 1/32$, $x_0 \in \sX$, $R>0$.
 It is sufficient to show that there exists $M$ (depending
only on $C_H$) such that if $B(z_i, 8 \delta R)$, $1\le i\le n$ are disjoint balls with centres
in $B(x_0,R)- B(x_0, R/4)$ then $n \le M$.
So let $B(z_i, 8 \delta R)$, $i=1, \dots ,n$ satisfy these conditions. 

Let $B_k = B(x_0, \half k \delta R)$ for $1/(2\delta) \le k \le 2/\delta$, and let 
$n_k$ be the number of balls $B(z_i, \delta R)$ which intersect
$\pd B_k$. Since each $B(z_i, \delta R)$ must intersect at least
one of the sets $\pd B_k$ we have $n \le \sum_k n_k$.
The previous Lemma gives $n_k \le C_1$, and thus $n \le 2 C_1/\delta$. 
\qed

\ms We now compare $g_D$ in two domains.

\begin{lemma} \label{L:domcap}
There exists a constant $C_0$ such that if $B=B(x_0,R)$ and
$2B=B(x_0,2R)$ then
\be
 g_{2B}(x,y) \le C_0 g_B(x,y) \,\, \hbox{ for } x,y \in B(x_0, R/4).
\ee
\end{lemma}

\proof Let $B'=B(x_0,R/2)$ and $y \in B(x_0,R/4)$. Choose
$x_1 \in \pd B'$ to maximise $g_{2B}(x_1,y)$. 
Let $\gam$ be a geodesic path from $x_0$ to $\pd B(x_0,3R)$, let $z_0$ be the point
on $\gam \cap \pd B$, and $A = B(z_0,R/4)$.

Using Lemma  \ref{L:hit} there exists $p_1>0$ such that
\begin{align*}
 p_A(w) = \bP^w( X_{\tau_B} \in A ) &\ge p_1, \q w \in B', \\ 
  \bP^z(  \tau_{2B} <  T_{B'}  ) &\ge p_1, \q z \in A. 
\end{align*}
Then
\begin{align*}
 g_{2B}(x_1,y) &= g_B(x_1,y) + \bE^{x_1} g_{2B}(X_{\tau_B},y) \\
  &= g_B(x_1,y) + \bE^{x_1} \one_{(X_{\tau_B} \in A)} g_{2B}(X_{\tau_B},y) 
  +  \bE^{x_1} \one_{(X_{\tau_B} \not \in A)} g_{2B}(X_{\tau_B},y) \\
 &\le   g_B(x_1,y) + p_A(x_1) \sup_{w \in A \cap \pd B} g_{2B}(w,y)
   + (1-p_A(x_1))  \sup_{z \in \pd B} g_{2B}(z,y) \\
&\le   g_B(x_1,y) + p_1 \sup_{w \in A \cap \pd B} g_{2B}(w,y)
   + (1-p_1)  \sup_{z \in \pd B} g_{2B}(z,y).
\end{align*}
If $w \in A$ then
$$ g_{2B}(w,y) = \bE^w \one_{( T_{B'} < \tau_{2B} )} g_{2B}( X_{T_{B'}},y)
 \le (1-p_1) \sup_{z \in \pd B'} g_{2B}(z,y) \le  (1-p_1) g_{2B}(x_1,y). $$
The maximum principle implies that $g_{2B}(z,y) \le g_{2B}(x_1,y)$
for all $z \in \pd B$.
Combining the inequalities above gives 
\begin{align*}
 g_{2B}(x_1,y) &\le g_B(x_1,y) + p_1 (1-p_1)  g_{2B}(x_1,y)
  + (1-p_1)  g_{2B}(x_1,y),
\end{align*}
which implies that
\be
 g_{2B}(x_1,y) \le  p_1^{-2}  g_B(x_1,y).
\ee
Now let $x \in B(x_0,R/4)$.  By Corollary \ref{C:hg3}
$$  g_{2B}(x_1,y) \le p_1^{-2}  g_B(x_1,y) \le C g_B(x,y). $$
Hence
\begin{align*}
 g_{2B}(x,y) ) &= g_{B'}(x,y) + \bE^x g_{2B}(X_{\tau_{B'}},y) \\
  &\le g_{B'}(x,y) + g_{2B}(x_1,y) \le (1+C)  g_B(x,y).
\end{align*}
\qed

\begin{cor} \label{C:cap12}
Let $A \ge 4$. There exists $C_0=C_0(C_H,A)$ such that for $x \in \sX$,
$r>0$,
\be
 \Cpc_{B(x,2Ar)}(B(x,r)) \le  \Cpc_{B(x,Ar)}(B(x,r)) \le C_0  \Cpc_{B(x,2Ar)}(B(x,r)) .
\ee
\end{cor}

\proof The first inequality is immediate from the monotonicity of capacity, and the
second one follows immediately from Lemmas \ref{L:cap1} and \ref{L:domcap}. \qed

\begin{lemma}
Let $A\ge 8$, and $D$ be a bounded domain in $\sX$. 
Let $x \in \sX$ and $r>0$ be such that $B(x,4r) \subset D$. \\
(a) There exists $C_0=C_0(C_H)$ such that
$$  \Cpc_{D}(B(x,r)) \le C_0  \Cpc_{D}(B(y,r)) \q \hbox{ for } y \in B(x,r). $$
(b) There exists $C_1=C_1(A,C_H)$ such that
\be
  \Cpc_{B(x,Ar)}(B(x,r)) \le C_1  \Cpc_{B(y,Ar)}(B(y,r))  \q \hbox{ for } y \in B(x,r). 
\ee
\end{lemma}

\proof (a) follows easily from Lemmas \ref{L:cap1} and \ref{L:domcap}. \\
(b) We have
\begin{align*}
   \Cpc_{B(x,Ar)}(B(x,r)) &\le C_2  \Cpc_{B(x,Ar)}(B(y,r))  \\
   &\le C_3  \Cpc_{B(x,2Ar)}(B(y,r)) \le C_1    \Cpc_{B(y,Ar)}(B(y,r)).
\end{align*}
\qed

We conclude this section with a capacity estimate which will
play a key role in our construction of a well behaved doubling measure.

\begin{prop} \label{P:capest}
Let $D \subset \sX$ be a bounded open domain, and
let $B(x_0,8R) \subset D$. Let $F \subset B(x_0,R)$. 
Let $b \ge 4$, and 
suppose there exist disjoint Borel sets $(Q_i, 1\le i \le n)$, with $n \ge 2$, such that
$$ F = \cup_{i=1}^n Q_i, $$
and for each $i$ there exists $z_i \in Q_i$ such that
$B(z_i, R/6b) \subset Q_i. $
Then there exists $\delta = \delta(b,C_H) >0$ such that
$$ \textstyle \Cpc_D(F) \le  (1-\delta) \sum_{i=1}^n  \textstyle \Cpc_D(Q_i).  $$
\end{prop}

\proof
Let $\nu_i$ and $h_i$ be the equilibrium measure and equilibrium 
potential respectively for $Q_i$, so that
and $h_i=  1$~q.e. on $Q_i$. Then
$$  \Cpc_D(Q_i) = \nu_i(\partial Q_i)= \nu_i(D). $$ 
By \eqref{e:cap-hit} and Lemma \ref{L:hit}, there exists $c_1>0$ such that
$$ h_i(y)  \ge c_1 \q \hbox{ for } y \in B(x_0,R) ~\mbox{q.e.} $$

Let $h =  \sum_{i=1}^n h_i$. Let $y \in F$, so that
there exists $i$ such that $y \in Q_i$. Then since $n \ge 2$, 
\bes
h(y)= \sum_{i=1}^n h_i(y)  
\ge 1  + \sum_{j\neq i} c_1 \ge 1+c_1,\q \hbox{ for } y \in F~\mbox{q.e.}
\ees
Consequently if $h' = \left[(1+c_1)^{-1}  h \right] \wedge 1$ then
$h' = 1$ quasi everywhere on $F$. It follows that
\begin{align*}
\textstyle \Cpc_D(F) &\le \sE(h',h') \le \sE(h',(1+c_1)^{-1}  h)
=(1+c_1)^{-1} \sum_{i=1}^n \int_{D} h' \,d\nu_i  \\
&\le  (1+c_1)^{-1}\sum_{i=1}^n \nu_i(D)= (1+c_1)^{-1}\sum_{i=1}^n  \Cpc_D(Q_i).
\end{align*}
The first inequality above follows from \eqref{e:cap-char}, the second inequality follows from the fact that $h'$ is a potential (see \cite[Corollary 2.2.2 and  Lemma 2.2.10]{FOT}), 
the third equality follows from \eqref{e:capmeas1} and 
the fourth inequality holds since $h' \le 1$.
\qed

\begin{remark}\label{r-local}
{\rm All the results in this section  can be localized in the 
following sense: if we assume the EHI holds at small scales (i.e. for radii less than 
some $R_1$)
then the conclusions of the results in this section also hold at 
sufficiently small scales.
}\end{remark}

\section{Construction of good doubling measures}
\label{sec:goodVD}

We continue to consider a length metric measure space with Dirichlet
form $(\sX,d,m,\sE,\sF^m)$ which satisfies the EHI and Assumptions \ref{A:green} and \ref{A:BG}.
The space $(\sX,d)$ satisfies metric doubling by Theorem \ref{T:BC}, 
and therefore by \cite{VK,LuS} there exists a doubling measure 
$\mu$ on $(\sX,d)$. However, this measure might be somewhat pathological 
(see  \cite[Theorem 2]{Wu}), and to prove the EHI 
we will require some additional regularity properties of $\mu$.
 In this section we adapt the argument of \cite{VK} to obtain a `good'
doubling measure, that is one which connects measures and capacities  of balls
in a satisfactory fashion.

\begin{definition}\label{d:goodmeas}
{\rm Let $D$ be either a ball $B(x_0,R) \subset \sX$ or the whole space $\sX$. If
$D=\sX$ fix $x_0 \in \sX$.
Let $C_0<\infty$ and $0< \beta_1 \le \beta_2$.
We say a measure $\nu$ on $D$ is {\em $(C_0, \beta_1, \beta_2)$-capacity good} if the following holds.
\begin{enumerate}[(a)]
\item The measure $\nu$ is doubling on all balls contained in $D$, that is 
\begin{equation} \label{e:m01}
	\frac{\nu\left(B(x,2s) \right)}{\nu\left(B(x,s)  \right)} \le C_0
	 \mbox{ whenever $B(x,2s) \subset D$. }
\end{equation}
\item 
For all $x \in D$ and $0<s_1<s_2$  such that $B(x,s_2)\subset D$,
\begin{equation} \label{e:m02}
 C_0^{-1} \left( \frac{s_2}{s_1}\right)^{\beta_1} \le
\frac{\nu(B(x,s_2)) \Cpc_{B(x,8s_1)} (B(x,s_1))} {\nu(B(x,s_1)) \Cpc_{B(x,8s_2)}(B(x,s_2))}
\le C_0 \left( \frac{s_2}{s_1}\right)^{\beta_2}.
\end{equation}
\item The measure $\nu$ is  absolutely continuous with respect to $m$  and we have
\begin{align}  \label{e:m03}
\esssup_{y \in B(x,1)} \frac{d \nu}{dm}(y) &\le C_0 \essinf_{y \in B(x,1)} \frac{d \nu}{dm}(y)
  \, \hbox{   whenever  $B(x,1) \subset D$}, \\
\label{e:m04}
	C_0^{-1-d(x_0,y)} &\le  \frac{d \nu}{dm}(y) \le C_0^{1+d(x_0,y)},
	\mbox{ for $m$-almost every $y \in D$.}
\end{align}
\end{enumerate}
}\end{definition}

The following  is the main result of this section.

\begin{theorem}[Construction of a doubling measure] \label{T:meas}
Let $(\sX,d)$ be a complete, locally compact, length metric space with a strongly local 
regular Dirichlet form $(\sE,\sF^m)$ on $L^2(\sX,m)$ which satisfies
Assumptions \ref{A:green} and \ref{A:BG} and the EHI.
Then there exist constants $C_0>1$,  $0< \beta_1 \le \beta_2$
and a measure $\mu$ on $\sX$ which is $(C_0, \beta_1, \beta_2)$-capacity good.
\end{theorem}

We begin by adapting the argument in \cite{VK} to construct
measure with the desired properties in  a ball $B(x_0,R)$. 
We then follow \cite{LuS} and 
obtain $\mu$ as a weak$^*$ limit of measures defined on an
increasing family of balls.

\begin{proposition}[Measure in a ball] \label{P:ballmeas}
Let $(\sX,d,m,\sE,\sF^m)$ be as in the previous Theorem.
There exist $C_0>1$,  $0< \beta_1 \le \beta_2$ such that for any ball
$B_0=B(x_0,r)  \subset \sX$ with $r\ge r_0$
there exists a measure
$\nu=\nu_{x_0,r}$ on $B_0$ which is 
$(C_0, \beta_1, \beta_2)$-capacity good.
\end{proposition}

The proof uses a family of generalized dyadic cubes, which 
provide a family of nested partitions of a space.

\begin{lemma}(\cite[Theorem 2.1]{KRS}) \label{L:gdc}
Let $(\sX,d)$ be a complete, length metric space satisfying \BC\hspace{0.1mm} 
and let $A \ge 4$ and $c_A=  \frac{1}{2} - \frac{1}{A-1}$.
Let $B_0= B(x_0,r)$ denote a closed ball in $(\sX,d)$. Then there exists a collection $\Sett{Q_{k,i}}{k \in \bZ_+, i \in I_k \subset \bZ_+}$  of Borel sets satisfying the following properties:
\begin{enumerate}[(a)]
\item  $B_0 = \cup_{i\in I_k} Q_{k,i}$ for all $k \in \bZ_+$.
\item If $m \le n$ and $i \in I_n$, $j \in I_m$ then either 
$Q_{n,i} \cap Q_{m,j} = \emptyset$ or else $Q_{n,i} \subset Q_{m,j}$.
\item For every $k  \in \bZ_+$, $i \in I_k$, there exists $x_{k,i}$ such that
\[
B(x_{k,i},c_A r A^{-k}) \cap B_0 \subset Q_{k,i} \subset B(x_{k,i},  A^{-k} r).
\]
\item The sets $N_k= \Sett{x_{k,i}}{i \in I_k}$, where $x_{k,i}$ are as defined in (c) 
are increasing, 
$N_0= \set{x_0}$, and $Q_{0,0}=B_0$. Moreover for each $k \in \bZ_+$
$N_k$ is a maximal $rA^{-k}$-separated subset ($rA^{-k}$-net) of $B_0$.
\item Property (b) defines a partial order $\prec$ on $\mathcal{I}=\Sett{(k,i)}{k \in \bZ_+, i \in I_k}$ by inclusion, where
$(k,i) \prec (m,j)$ whenever  $Q_{k,i} \subset Q_{m,j}$.
\item There exists $C_M>0$ such that, for all $k \in \bZ_+$ and 
for all $x_{k,i} \in N_k$, the  `successors'
\[
S_{k}(x_{k,i})= \Sett{x_{k+1,j}}{(k+1,j) \prec (k,i)}
\]
satisfy
\be \label{e:in1}
C_M \ge \abs{S_k(x_{k,i})}  \ge 2.
\ee
Moreover, by property (c), we have $d(x_{k,i},y) \le A^{-k}r$ for all $y \in S_k(x_{k,i})$.
\end{enumerate}
\end{lemma}

\sms
We now set $A=8$ and 
until the end of the proof of Proposition \ref{P:ballmeas} 
we fix a ball $B_0= B(x_0,r)$.
We remark that the constants in the rest of the section do not depend on the ball 
$B_0$: they depend only on the constants in EHI and \BC.

We fix a family 
\[ \Sett{Q_{k,i}}{k \in \bZ_+, i \in I_k \subset \bZ_+},\] 
of generalized dyadic cubes as given by Lemma \ref{L:gdc},
and define the nets $N_k$ and successors $S_k(x)$ as in the lemma.
For $k\ge 1$, we define the predecessor $P_k(x)$ of $x\in N_k$ to be the unique 
element of $N_{k-1}$ such that $x \in S_{k-1}(P_k(x))$. 
Note that for $x \in N_k$, $S_k(x) \subset N_{k+1}$ whereas $P_k(x) \in N_{k-1}$.  For $x \in B_0$, we denote by $Q_k(x)$ the unique  $Q_{k,i}$ such that $x \in Q_{k,i}$.
 For $x \in N_k$, we denote by $c_k$ the capacity
\[
c_k(x)= \operatorname{Cap}_{B(x, A^{-k+1}r)} (Q_{k}(x)).
\]
The following lemma provides useful estimates on $c_k$.

 \begin{lemma}[Capacity estimates for generalized dyadic cubes] \label{l-ced}
 There exists $C_1>1$ such that the following hold.\\
(a)  For all $k \in \bZ_+$ and for all $x,y \in N_k$, 
such that $d(x,y) \le 4 rA^{-k}$, we have
\begin{equation}\label{e:ce1}
C_1^{-1} c_k(y) \le c_k(x) \le C_1 c_k(y).
\end{equation}
(b) For all $k \in \bZ_+$, for all $x \in N_k$, for all 
$y \in S_k(x)$, we have
\begin{equation}\label{e:ce2}
C_1^{-1} c_k(x) \le c_{k+1}(y) \le C_1 c_k(x).
\end{equation}
\end{lemma}

\proof
	First, we observe that there is $C>1$ such that 
\begin{equation} \label{e:ce3}
	 C^{-1} \left( g_{B(x,A^{-k+1}r)}(x, A^{-k}r)\right)^{-1} 
	 \le c_k(x) \le C \left(g_{B(x,A^{-k+1}r)}(x, A^{-k}r)\right)^{-1}
\end{equation}
	for all $x \in B(x_0,r)$.
The upper bound in \eqref{e:ce3} follows from Lemma \ref{L:gdc}(c), domain monotonicity of capacity and  Lemma \ref{L:cap1}.
For the lower bound, we again use Lemma \ref{L:gdc}(c) to choose a point 
$z \in \gamma(x_0,x) \cap B_0$ such that $d(x,z)=c r A^{-k}/2$  where $c$ is as given by Lemma \ref{L:gdc}(c). 
By the triangle inequality $Q_k(x) \supset B(z,c r A^{-k}/2)$. 
The lower bound again follows from domain monotonicity, Lemma \ref{L:cap1} and standard chaining arguments using EHI.  
The estimates \eqref{e:ce1} and \eqref{e:ce2} then follow from \eqref{e:ce3}, 
domain monotonicity of capacity and Lemma \ref{L:domcap}.
\qed

We record one more estimate regarding the subadditivity of $c_k$,
which will play an essential role in  ensuring \eqref{e:m02}.

\begin{lemma}[Enhanced subaddivity estimate] \label{L:ESP}
There exists $\delta \in (0,1)$ such that for all $k \in \bZ_+$, for all $x \in N_k$, we have
\[
c_k(x) \le (1-\delta) \sum_{y \in S_k(x)} c_{k+1}(y).
\]
\end{lemma}
\proof
By the triangle inequality, $B(y,A^{-k}r) \subset B(x,A^{-k+1}r)$ 
for all $k \in \bZ_+, x \in N_k, y \in S_k(x)$.
The lemma now follows from Proposition \ref{P:capest}
and domain monotonicity of capacity.
\qed

We now follow the Vol'berg-Konyagin construction closely, but with some essential changes. 
Recall that we want to construct a doubling measure $\mu$ on $B_0$ satisfying the 
estimates in Definition \ref{d:goodmeas}.

\begin{lemma}\label{l-ind} (See \cite[Lemma, p. 631]{VK}.)
Let $B_0= B(x_0,r)$ and let $c_k$ denote the capacities of the corresponding generalized dyadic cubes as defined above.
There exists $C_2 \ge 1$ satisfying the following. 
Let $\mu_k$ be a probability measure on $N_k$ such that
\be \label{e:muk}
\frac{\mu_k(e')}{c_k(e')} \le  C_2^2 \frac{\mu_k(e'')}{c_k(e'')}
 \q \hbox{for all $e',e'' \in N_k$ with $d(e',e'') \le 4 A^{-k}r$}.
\ee
Then there exists a 
probability measure $\mu_{k+1}$ on $N_{k+1}$ such that
\begin{enumerate}[(1)]
\item 
For all $g',g'' \in N_{k+1}$ with $d(g',g'') \le 4 A^{-k-1}r$ we have
\be  \label{e:mure1}
\frac{\mu_{k+1}(g')}{c_{k+1}(g')} \le  C_2^2 \frac{\mu_{k+1}(g'')}{c_{k+1}(g'')}.
\ee
\item 
Let $\delta \in (0,1)$ be the constant in Lemma \ref{L:ESP}.
For all points $e \in N_k$ and $g \in S_k(e)$,
\be \label{e:mure2}
C_2^{-1}\frac{\mu_{k}(e)}{c_{k}(e)} \le    \frac{\mu_{k+1}(g)}{c_{k+1}(g)} 
\le (1-\delta) \frac{\mu_{k}(e)}{c_{k}(e)}.
\ee
\item The construction of the measure $\mu_{k+1}$ from the measure 
$\mu_k$ can be regarded as the transfer of masses from the points 
$N_k$ to those of $N_{k+1}$, with no mass transferred over a distance greater 
than $(1+4/A)A^{-k}r$.
\end{enumerate}
\end{lemma}

\begin{remark} {\rm
The key differences from the Lemma in \cite{VK} are, first, that we require 
the relations \eqref{e:muk}, \eqref{e:mure1} and \eqref{e:mure2} for the 
ratios $\mu_k/c_k$ rather than just for $\mu_k$, and second, the presence
of the term $1-\delta$ in the right hand inequality in \eqref{e:mure2}.
}\end{remark}

\proof
By Lemma \ref{L:gdc}(f) we have $|S_k(x)| \le C_M$ for all $x$, $k$.
Set
\[
C_2=C_1 C_M,
\]
where $C_1$ is the constant in \eqref{e:ce1}.
Let $\mu_k$ be a probability measure $N_k$ satisfying \eqref{e:muk}.

Let $e \in N_k$; we will construct 
$\mu_{k+1}(g)$ for $g \in S_k(e)$ by mass transfer. 
Initially we distribute the mass $\mu_k(e)$ to $g \in S_k(e)$
so that the mass of $g\in S_{k+1}(e)$ is proportional to $c_{k+1}(g)$.
We therefore set
\[
f_0(g)= \frac{c_{k+1}(g)}{\sum_{g' \in S_k(e)} c_{k+1}(g')} \mu_k(e), 
\q \mbox{ for all $e \in N_k$ and $g \in S_k(e)$.}
\]
By \eqref{e:in1}, Lemma \ref{l-ced} and Lemma \ref{L:ESP}, we have
 \begin{equation} \label{e:in2}
C_2^{-1}\frac{\mu_{k}(e)}{c_{k}(e)} \le    \frac{f_0(g)}{c_{k+1}(g)} 
\le (1- \delta) \frac{\mu_{k}(e)}{c_{k}(e)},
\end{equation}
for all points $e \in N_k$ and $g \in S_k(e)$.
If the measure $f_0$ on $N_{k+1}$ satisfies condition (1) of the Lemma, 
we set $\mu_{k+1}=f_0$. Condition (2) is satisfied by \eqref{e:in2}, 
and (3) is obviously satisfied by Lemma \ref{L:gdc}(c).

If $f_0$ does not satisfy condition (1) of the Lemma, 
then we proceed to adjust the masses of the points in $N_{k+1}$ in the
following fashion.
Let $p_1,\ldots,p_T$ be the pairs of points
$\set{g',g''}$ with $g',g'' \in N_{k+1}$ with $0 < d(g',g'') \le 4 A^{-k-1}r$. 
We begin with the pair $p_1= \set{g_1',g_1''}$. If
$$ \frac{f_0(g_1')}{ c_{k+1}(g_1')}\le C_2^2 \frac{f_0(g_1'')}{ c_{k+1}(g_1'')},
\q \hbox{ and } 
\frac{f_0(g_1'')}{ c_{k+1}(g_1'')} \le C_2^2 \frac{f_0(g_1')}{ c_{k+1}(g_1')}, $$
then we set $f_1=f_0$. If one of the inequalities is violated, say the first,
then we define the measure $f_1$ by a suitable transfer of mass
from $g_1'$ to $g_1''$. 
We set $f_1(g) = f_0(g)$ for $g \neq g_1',g_1''$, and set
$f_1(g_1') =f_0(g_1') - \alpha_1$,  $f_1(g_1'') = f_0(g_1'') + \alpha_1$,
where $\alpha_1 >0$ is chosen such that
$$ \frac{f_1(g_1')}{ c_{k+1}(g_1')} = C_2^2 \frac{f_1(g_1'')}{ c_{k+1}(g_1'')}.$$
We then consider the pair $p_2$, and construct the measure $f_2$ from $f_1$ in exactly the same way,
by a suitable mass transfer between the points in the pair if this is necessary.
Continuing we obtain a sequence of measures $f_j$, and we find that 
$\mu_{k+1}=f_T$ is the desired measure in the lemma.

The proof that $\mu_{k+1}$ satisfies the properties (1)--(3) is almost the same as in \cite{VK}.
We note that a key property of the construction is that we cannot have chains of mass transfers:
as in  \cite{VK} there are no pairs $p_j=(g_1, g_2)$, $p_{j+i}=(g_2,g_3)$
such that at step $j$ mass is transferred from $g_1$ to $g_2$, and then at a 
later step $j+i$ mass is transferred from $g_2$ to $g_3$. (See \cite[p. 633]{VK}.)
\qed

To construct the doubling measure in Proposition \ref{P:ballmeas} 
we use Lemma \ref{l-ind} for large scales, and rely on \BG\hspace{0.1mm} for small scales. 

{\sm {\em Proof of Proposition \ref{P:ballmeas}.}} 
Recall that  $A=8$.
Let $\mu_0$ be the probability measure on $N_0=\set{x_0}$. We use Lemma \ref{l-ind} 
to inductively construct probability measures $\mu_k$ on $N_k$.
For $x \in B_0$, by $E_k(x)$ we denote the unique 
$y \in N_k$ such that $Q_k(x)=Q_k(y)$. Note that by 
the construction
\begin{equation}\label{e:ms1}
d(x,E_k(x))< A^{-k}x, \hspace{4mm}  P_{k+1}(E_{k+1}(x))= E_k(x),
\end{equation}
for all $x \in B_0$ and for all $k \in \bZ_+$. 
Let $l$ denote the smallest non-negative integer such that $A^{-l}r \le r_0/A^2$;
since $r \ge r_0$ we have $l \ge 2$.
The desired measure $\nu=\nu_{x_0,r}$ is given by
\[
f(z)= \alpha \sum_{y \in N_l} \frac{\mu_l(y)}{m(Q_l(y))}\one_{Q_l(y)}(z), 
\q \nu(dz)= f(z) \,m(dz),
\]
where $\alpha >0$ is chosen so that $f(x_0)=1$.
Note that we have $\mu_l(x)= \alpha^{-1} \nu(Q_l(x))$ for all $x \in N_l$.

First, we show \eqref{e:m03} and \eqref{e:m04}. By the argument in  \cite[Lemma 2.5]{Kan} 
there exists $C_3>1$ (that does not depend on $r$) such for any pair of points 
$x,y \in B_0$ that can be connected by a geodesic that stays within $B_0$, 
there exists sequence of points 
$E_l(x)=y_0,y_1,\ldots,y_{n-1},y_n=E_l(y)$ in $N_l$, with 
$n \le C_3 (1+d(x,y))$ and $d(y_i,y_{i+1}) \le 4 A^{-l} r$ for all $i=0,1,\ldots,n-1$. 
By comparing successive $\mu_l(y_i)$ using Lemma \ref{l-ind}(3) and by 
comparing successive $m(Q_l(y))$ using the volume doubling property of $m$  
at small scales \eqref{e:vds}, we obtain \eqref{e:m03} and \eqref{e:m04}.

For  rest of the proof we can 
without loss of generality assume that $\alpha=1$ in the definition of $\nu$.
If $x\in N_k$ then in the mass transfer from $\mu_k$ to $\mu_l$ each
piece of mass moves a distance at most
$ (1 + 4A^{-1}) \sum_{i=k}^l A^{-i} r $. An additional 
distance of at most $A^{-l}r$ is then travelled in the transfer from
$\mu_l$ to $\nu$. Since $A \ge 8$ 
the mass $\mu_k(x)$ from $x \in N_k$ travels a distance of at most 
\begin{equation} \label{e:mtran}
  (1 + 4A^{-1}) \sum_{i=k}^l A^{-i} r + A^{-l}r  < 2 A^{-k} r.
\end{equation}

Next we show that there exists $C_4$ such that
\begin{equation}\label{e:ms2}
 \mu_{M+1} \left(E_{M+1}(x)\right) \le \nu(B(x,s)) \le C_4  \mu_{M+1}(E_{M+1}(x)).
\end{equation}
for all $x \in B_0$ and for all $A^{-l+1}r <s<r$. Here $M = M(s) \in \bZ_+$ 
is the unique integer such 
that $s/A \le A^{-M}r < s$. Note that $M \le l-1$.

By \eqref{e:mtran} the mass transfer of 
$\mu_{M+1} \left(E_{M+1}(x)\right)$ from the point $E_{M+1}(x)$ 
takes place over a distance at most
$ 2 A^{-M-1} r \le {2s}/{8}$.
Since $d(x,E_{M+1}(x)) \le A^{-M-1}r < s/8$, the triangle inequality gives the lower bound
in \eqref{e:ms2}.

To prove the upper bound, 
recall from \eqref{e:mtran} that none of the mass in $N_{M-1} \setminus B(x, s+ 2 A^{-M+1}r)$ of $\mu_{M-1}$ falls in $B(x,s)$. This implies that
\begin{equation}\label{e:ms4}
\nu(B(x,s)) \le \mu_{M-1} \left( N_{M-1} \cap  B(x, s+ 2 A^{-M+1}r) \right).
\end{equation}
Since $s \le A^{-M+1}r$ and $N_{M-1}$ is an $A^{-M+1}r$-net, by \BC\hspace{0.1mm}
 there exists $C_5>1$ such that
\begin{equation} \label{e:ms5}
\abs{ N_{M-1} \cap  B(x, s+ 2 A^{-M+1}r)} \le C_5.
\end{equation}
By the triangle inequality $d(x,E_{M-1}(x),y) < 4 A^{-M+1}r$ for all 
$y \in   N_{M-1} \cap  B(x, s+ 2 A^{-M+1}r)$. Therefore by \eqref{e:ms4}, \eqref{e:ms5}, Lemma \ref{l-ind}, Lemma \ref{l-ced}, there exists $C_6 >1$ such that
 \begin{equation}\label{e:ms6}
\nu(B(x,s)) \le C_6 \mu_{M-1}(E_{M-1}(x)).
\end{equation}
Combining \eqref{e:ms6} along with \eqref{e:ms1} and Lemma \ref{l-ced}, 
we obtain the desired upper bound in \eqref{e:ms2}. 

For small balls we rely on \BG\hspace{0.1mm} as follows. 
If $B(x,s) \subset B_0$, $s \le A^{-l+2}r$, $y \in B(x,s)$ there exists $C_7>1$ 
such that $E_l(x)$ and $E_l(y)$ can be connected by a chain of points in $N_l$ 
given by $E_l(x)=z_0,z_l,\ldots,z_{N-1},z_N=E_l(y)$ with $N \le C_7$. 
This can be shown by essentially the same argument as  \cite[Lemma 2.5]{Kan} 
or \cite[Proposition 2.16(d)]{MS1}.  Combining this with \eqref{e:vds},
Lemmas \ref{l-ind} and \ref{l-ced} 
we obtain that there exists $C_8>1$ such that
\[
C_8^{-1} f(x) \le f(y) \le C_8 f(x).
\]
Therefore for all balls $B(x,s) \subset B_0$ with  $s < A^{-l+2} r$, we have
\begin{equation}\label{e:ms7}
 C_8^{-1} f(x) m(B(x,s)) \le \nu(B(x,s)) \le C_8 f(x) m(B(x,s)).
\end{equation}

Combining \eqref{e:ms2} with Lemmas \ref{l-ind} and \ref{l-ced}, we obtain the volume doubling property 
for $\nu$ for all balls whose radius $s$ satisfies $A^{-l+1}r<s<r$. 
The estimate \eqref{e:ms7} and \BG\hspace{0.1mm} for the measure $m$ implies 
the volume doubling property for balls $B(x,s)$ with $s\le A^{-l+1}r$ and 
$B(x,2s) \subset B_0$.  This completes the proof of the doubling property 
given in \eqref{e:m01}.

It remains to verify \eqref{e:m02}. 
By an application of EHI, \eqref{e:ce2}, \eqref{e:ce3} 
along with Lemmas \ref{L:cap1}, \ref{L:domcap} and  \ref{L:greenc1}  there exists $C_9>1$ such that 
\begin{equation}\label{e:ms8}
C_9^{-1} c_{M+1}(E_{M+1}(x)) 
 \le  \Cpc_{B(x,As)}\left(B(x,s)\right) \le C_9 c_{M+1}(E_{M+1}(x))
\end{equation} 
for all $x \in B_0$, for all $A^{-l+1}r < s \le r$, where $M=M(s)$ is as above.

The equations \eqref{e:ms2}, \eqref{e:ms7} and \eqref{e:ms8} link
the $\nu$-measure and capacity of balls with those of the generalized cubes
$Q_{k,i}$, while Lemma \ref{l-ced} and Lemma \ref{l-ind}
link $\nu$-measures and capacities of the $Q_{k,i}$ with their
successors and neighbours. Using these links, as well as the regularity
on small scales given by Assumption \ref{A:BG}, \eqref{e:m02}
follows by a straightforward argument. \qed

{\sm {\em Proof of Theorem \ref{T:meas}.}} 
Fix $x_0 \in \sX$.
For $n\ge 1 \vee r_0$ let $\nu_{x_0,n}$ be the measure given by 
Proposition \ref{P:ballmeas}, and let
\[
f_n:=\frac{d\nu_{x_0,n}}{dm}.
\]
Then $f_n \in L^\infty(B(x_0,n), m)$ and by \eqref{e:m03} we have
for $1 \le k \le n$ that
\begin{equation}\label{e:tms1}
C_0^{-1-k} \le \essinf_{B(x_0,k)} f_n \le \esssup_{B(x_0,k)} f_n \le C_0^{1+k}.
\end{equation}
A compactness argument similar to that in  \cite{LuS} yields the existence
of a subsequence ${n_k}$ and a measurable function $f$, bounded on
compacts, such that
\begin{equation} \label{e:tms2}
\int_{\sX} h f \, dm = \lim_{k \to \infty}\int_{\sX} h f_{n_k} \, dm
\end{equation}
for all $h \in L^1(\sX,m)$ with compact support.
Taking $d \mu = fdm$ then yields the required measure. \qed

\section{ Quasisymmetry, time change, and the EHI} \label{sec:qs}

We recall the definition of quasisymmetry -- see \cite{Hei} for a nice exposition to this theory. 
For many results in this section, we do not require that the metric space $(\sX,d)$ is a
length space.

\begin{definition}\label{d:qs}
{\rm	A \emph{distortion function} is a homeomorphism of $[0,\infty)$ onto itself. 
	Let $\eta$ be a distortion function. 
	A map $f:(\sX_1,d_1) \to (\sX_2,d_2)$ between metric spaces is said to be
	 \emph{$\eta$-quasisymmetric}, if $f$ is a homeomorphism and
	\[
 \frac{d_2(f(x),f(a))}{d_2(f(x),f(b))} \le \eta\left(\frac{d_1(x,a)}{d_1(x,b)}\right)  
	\]
	 for all triples of points $x,a,b \in \sX_1, x \neq b$. We say $f$ is a \emph{quasisymmety} if it is $\eta$-quasisymmetric for some distortion function $\eta$.
	 We say that metric spaces $(\sX_1,d_1)$ and $(\sX_2,d_2)$ are quasisymmetric, if there exists a quasisymmetry $f:(\sX_1,d_1) \to (\sX_2,d_2)$.
	 We say that  metrics $d_1$ and $d_2$  on $\sX$ are \emph{quasisymmetric}, if the identity map $\operatorname{Id}:(\sX,d_1) \to (\sX,d_2)$ is a quasisymmety. }
\end{definition}

If $f:(\sX_1,d_1) \to (\sX_2,d_2)$ is $\eta$-quasisymmetric, then $f^{-1}:(\sX_2,d_2) \to (\sX_1,d_1)$ is $\zeta$-quasisymmetric, where $\zeta(t)=1/\eta^{-1}(1/t)$.
Quasisymmetry is an equivalence relation among metric spaces \cite[Proposition 10.6]{Hei}.
The following comparison of annuli follows readily from the definition.
\begin{lem} (See \cite[Lemma 1.2.18]{MT}) \label{l:cann}
Let the identity map $\operatorname{Id}:(\sX,d_1) \to (\sX,d_2)$ be an $\eta$-quasisymmetry 
for some distortion function $\eta$. 
Then for all $A>1,x \in \sX, r>0$, there exists $s>0$ such that, writing $B_i$ for balls in $(\sX,d_i)$
\be \label{e:ann1}
B_2(x,s) \subset B_1(x,r) \subset B_1(x,Ar) \subset B_2(x, \eta(A)s).
\ee
Moreover, for all $A >1,x \in \sX,r >0$, there exists $s>0$ such that
\be \label{e:ann2}
B_1(x,r) \subset B_2(x,s) \subset B_2(x,As) \subset B_1(x,A_1r),
\ee
where $A_1= {1}/{\eta^{-1}(A^{-1})}$.
\end{lem}

The property of being a length metric space is not preserved under 
a quasisymmetric change of metric.
Nevertheless, many properties that are relevant to heat kernel estimates and Harnack inequalities
 are preserved under such transformations. 
For instance, it is well known that the metric doubling property (MD) is  a quasisymmetry invariant \cite[Theorem 10.18]{Hei}.
It follows easily from Lemma \ref{l:cann} that
 quasisymmetric metrics have the same doubling measures. 
 The EHI is also a quasisymmetry invariant as shown in the following easy but important lemma.

\begin{lem}  \label{l:ehiqs}
	Let $(\sX,d_i,\mu,\sE,\sF^\mu), i=1,2$ be two MMD spaces such that $d_1$ and $d_2$ are quasisymmetric. If $(\sX,d_2,\mu,\sE,\sF^\mu)$ satisfies EHI, then so does $(\sX,d_1,\mu,\sE,\sF^\mu)$.
\end{lem}
\proof
Let $C_H,A>1$ be constants corresponding to EHI for $(\sX,d_2,\mu,\sE,\sF^\mu)$.
Then by \eqref{e:ann2}, we have EHI for $(\sX,d_2,\mu,\sE,\sF^\mu)$ with constants $C_H,A_1>1$, 
where $A_1$ is  as given in Lemma \ref{l:cann}. \qed

We introduce the notion of a regular scale function.
\begin{definition} \label{d:regscale} {\rm
	We say that a function $\Psi:\sX \times [0,\infty) \to [0,\infty)$ on a metric space $(\sX,d)$ is a \emph{regular scale function} if $\Psi(x,0)=0$ for all $x$ and 
	there exist $C_1,\beta_1,\beta_2>0$ such that, for all $x,y \in \sX$,  $0<s \le r$ we have,
	writing $d(x,y)=R$, 
	\be \label{e:Psireg1}
	C_1^{-1} \Big(   \frac{r}{R \vee r}   \Big)^{\beta_2}  \Big( \frac{R \vee r}{s} \Big)^{\beta_1}
	\le  \frac{ \Psi(x,r)}{\Psi(y,s)} 
	\le C_1 \Big(   \frac{r}{R \vee r}   \Big)^{\beta_1}  \Big( \frac{R \vee r}{s} \Big)^{\beta_2}.  
	\ee 
}\end{definition}

We now recall the notion of \emph{uniform perfectness} \cite[Section 11.1]{Hei}.

\begin{definition} \label{d:up} {\rm
(1) A metric space $(\sX,d)$ is \emph{uniformly perfect} if there exists $C > 1$ so that for each $x \in \sX$, 
and for each $r>0$, the set $B(x,r)\setminus B(x,r/C)$ is nonempty whenever $\sX \setminus B(x,r)$ is non-empty. \\
(2) 
A measure $\mu$  satisfies the {\em reverse doubling property} (RVD) if there exist
constants $C_0$ and $\alpha>0$ such that
\be
  \frac{ \mu(B(x,r))} { \mu(B(x,s))} \ge C_0 (r/s)^\alpha \q \hbox{ for } x \in \sX, 0<s \le r.
  \ee
}\end{definition}

\begin{remark} \label{R:up}
{\rm Every connected metric space is uniformly perfect.  
Uniform perfectness is a quasisymmetry invariant -- see \cite[Exercise 11.2]{Hei}. 
If $\mu$ satisfies (VD) and $(\sX,d)$ is uniformly perfect, then $\mu$ satisfies (RVD) -- see  \cite[Exercise 13.1]{Hei}.
}\end{remark}

Next, we associate a quasisymmetric metric $d_\Psi$ to any regular scale function 
$\Psi$ on $(\sX,d)$, such that $d_\Psi$ relates nicely to $\Psi$.

\begin{proposition} \label{p:metric}
	Let $\Psi$  be a regular scale function on a metric space $(\sX,d)$. 
There exists a metric $d_\Psi:\sX \times \sX \to [0,\infty)$ satisfying the following properties:
\begin{enumerate}[(a)]
	\item There exist $C,\beta >0$ such that, for all $x , y \in \sX$ we have
	\be \label{e:defbeta}
	C^{-1} \Psi(x,d(x,y)) \le	d_\Psi(x,y)^\beta \le C \Psi(x,d(x,y)).
	\ee
	\item $d$ and  $d_\Psi$ are quasisymmetric.
	\item Assume in addition that $(\sX,d)$  (or equivalently $(\sX,d_\Psi)$) has infinite diameter and is uniformly perfect. Fix $A>1$.  Let $B_\Psi$ and $B$ denote metric balls in $(\sX,d_\Psi)$ and $(\sX,d)$ respectively. If either $B_\Psi(x,s) \subset B(x,r) \subset B_\Psi(x,As)$ or $B(x,r) \subset B_\Psi(x,s) \subset B(x,Ar)$ holds for some $x \in \sX, r>0, s>0$, then there is a constant $C_1>1$ (which does not depend on  $x \in \sX, r>0, s>0$)  such that
\be \label{e:qsm}
C_1^{-1} s^\beta \le \Psi(x,r) \le C_1 s^\beta,
\ee
where $\beta >0$ is as given by \eqref{e:defbeta}.
\end{enumerate}
\end{proposition}
\proof
(a)
Let $D(x,y) = \Psi(x,d(x,y))+ \Psi(y,d(x,y))$. Using \eqref{e:Psireg1}
it is straightforward to check that
 there exists $K \ge 1$ such that $D(x,y) \le K(D(x,z)+D(z,y))$ for all $x,y,z \in \sX$.
Therefore by \cite[Proposition 14.5]{Hei} and \eqref{e:Psireg1}, there exists a metric 
$d_\Psi$ on $\sX$ and $\beta>0$ that satisfies \eqref{e:defbeta}.\\
(b)
By \eqref{e:Psireg1} and \eqref{e:defbeta}, there exists $C_1>0$, $0<\gamma_1 \le \gamma_2$ such that
\[
C_1^{-1}\left(\frac{d(x,a)}{d(x,b)} \right)^{\gamma_1} \le \frac{d_\Psi(x,a)}{d_\Psi(x,b)} \le C_1 \left(\frac{d(x,a)}{d(x,b)} \right)^{\gamma_2}
\]
for all $x,a,b \in \sX$ that satisfy $d(x,a) \ge d(x,b)>0$. 
(We can take $\gam_i = \beta_i/\beta$.)
Therefore the identity map $\operatorname{Id}: (\sX,d) \to (\sX,d_\Psi)$ is quasisymmetric where the homeomorphism $\eta$ in Definition \ref{d:qs} can be chosen as $\eta(t)= C_1 \max(t^{\gamma_1},t^{\gamma_2})$.\\
(c) As the two cases are very similar, we just treat
the case $B_\Psi(x,s) \subset B(x,r) \subset B_\Psi(x,As)$.
By uniform perfectness, there exists $C_2>1$ such that 
there are points $y_1 \in B_\Psi(x,s) \setminus B_\Psi(x,s/C_2)$ and 
$y_2 \in B_\Psi(x,C_2As) \setminus B_\Psi(x,s)$. The upper bound in \eqref{e:qsm} follows from 
\[
\Psi(x,r) \le c \Psi(x,d(x,y_2)) \le c' d_\Psi(x,y_2)^\beta \le c'' s^\beta,
\]
where we used \eqref{e:Psireg1} and $d(x,y_2) \ge r$ in the first estimate, \eqref{e:defbeta} in the second, and $y_2 \in  B_\Psi(x,C_2As)$ in the final estimate.
Similarly, the lower bound in \eqref{e:qsm} follows from 
\[
\Psi(x,r) \ge c \Psi(x,d(x,y_1)) \ge c'  d_\Psi(x,y_1)^\beta \ge c''   s^\beta,
\]
where we used \eqref{e:Psireg1} and $d(x,y_1) \le r$ in the first estimate, \eqref{e:defbeta} in the second, and $y_1 \notin  B_\Psi(x,s/C_2)$ in the final estimate.
\qed 

We now introduce Poincar\'e and cutoff energy inequalities with respect to a regular scale function $\Psi$.
Recall that a \emph{cutoff function} $\vp$
for $B_1 \subset B_2$ is any function $\vp \in \sF^\mu$ such that $0\le \vp \le 1$ in $\sX$, $\vp \equiv 1$ 
in an open neighbourhood of $\overline{B_1}$, and $\operatorname{supp} \vp\Subset B_2$.

\begin{definition} \label{d:pi-cs}
	Let $(\sX,d,\mu,\sE,\sF^\mu)$ be a MMD space and let $\Psi$ be a regular scale function.
	
	 We say that $(\sX,d,\mu,\sE,\sF^\mu)$ satisfies the \emph{Poincar\'e inequality} \hypertarget{pi}{$\operatorname{PI}(\Psi)$}, if there exists constants $C,A \ge 1$ such that 
	for all $x\in \sX$, $R>0$ and $f \in \sF^\mu$
	\be  \tag*{$\operatorname{PI}(\Psi)$}
	\int_{B(x,R)} (f - \ol f)^2 \,d\mu  \le C \Psi(x,R) \, \int_{B(x,AR)}d\Gamma(f,f),
	\ee
	where $\ol f= \mu(B(x,r))^{-1} \int_{B(x,R)} f\, d\mu$.
	
We say that $(\sX,d,\mu,\sE,\sF^\mu)$ satisfies the 
\emph{cutoff energy inequality} \hypertarget{cs}{$\operatorname{CS}(\Psi)$}, 
if there exist $C_1,C_2>0,A>1$ such that the following holds.
		For all $R>0$, $x \in \sX$ with $B_1=B(x,R)$, $B_2=B(x,AR)$, there exists a cutoff function $\vp$
		for $B_1 \subset B_2$ such that for any $u \in \sF^\mu \cap L^\infty$,
		\be  \tag*{$\operatorname{CS}(\Psi)$}
		\int_{B_2 \setminus B_1} u^2 d\Gam(\vp,\vp) \le C_1 \int_{B_2 \setminus B_1} d\Gam(u,u)
		+ \frac{C_2}{\Psi(x,R)} \int_{B_2 \setminus B_1} u^2 \,d\mu.
		\ee 
	If  there exists $\beta >0$ such that $\Psi(x,r)=r^\beta$ for all $x \in \sX,r>0$, we denote the properties 	$\operatorname{PI}(\Psi)$ and $\operatorname{CS}(\Psi)$ by \hypertarget{pib}{$\operatorname{PI}(\beta)$} and \hypertarget{csb}{$\operatorname{CS}(\beta)$} respectively.
\end{definition}

The following lemma shows that the Poincar\'e and cutoff energy inequalities take a much simpler form 
with respect to the  metric $d_\Psi$. 
 
 \begin{lem} \label{l:picsqs}
	Let $(\sX,d,\mu,\sE,\sF^\mu)$ be an unbounded, uniformly perfect MMD space and let $\Psi$ be a regular scale function. Let $d_\Psi$ be the metric constructed in Proposition \ref{p:metric} with $\beta>0$ as given in \eqref{e:defbeta}. Then
	\begin{enumerate}[(a)]
		\item $(\sX,d,\mu,\sE,\sF^\mu)$ satisfies \hyperlink{pi}{$\operatorname{PI}(\Psi)$} if and only if $(\sX,d_\Psi,\mu,\sE,\sF^\mu)$ satisfies \hyperlink{pib}{$\operatorname{PI}(\beta)$}.
		\item $(\sX,d,\mu,\sE,\sF^\mu)$ satisfies \hyperlink{cs}{$\operatorname{CS}(\Psi)$} if and only if $(\sX,d_\Psi,\mu,\sE,\sF^\mu)$ satisfies \hyperlink{csb}{$\operatorname{CS}(\beta)$}.
	\end{enumerate}
\end{lem}
\proof
As before, we denote balls in the $d_\Psi$ and $d$ metrics by $B_\Psi$ and $B$ respectively.\\
(a) Let $(\sX,d,\mu,\sE,\sF^\mu)$ satisfy \hyperlink{pi}{$\operatorname{PI}(\Psi)$} with constants $C,A \ge 1$. By \eqref{e:ann2}, there exists $A'>1$ such that for all $x \in \sX, r>0$, there exists
$r'=r'(r)>0$ such that 
\be \label{e:pcq1}
B_\Psi(x,r) \subset B(x,r') \subset B(x,Ar') \subset B_\Psi(x,A'r).
\ee
Let $x \in X, r>0$ be arbitrary and let $r'>0, A'>1$ be given as above.
By \hyperlink{pi}{$\operatorname{PI}(\Psi)$}, \eqref{e:pcq1}, and Proposition \ref{p:metric}(c), there exists $C'>1$ such that
\be \label{e:pcq2}
\int_{B(x,r')} \abs{f-f_{B(x,r')}}^2\,d\mu \le  C \Psi(x,r') \int_{B(x,Ar')} d\Gamma(f,f) \le C' r^\beta \int_{B_\Psi(x,A'r)} d\Gamma(f,f),
\ee
for all $f \in \sF^\mu$.
Further, for all $f \in L^2(\sX,\mu)$
\begin{align} \label{e:pcq3}
\int_{B(x,r')} \abs{f-f_{B(x,r')}}^2\,d\mu &= \min_{a \in \bR}\int_{B(x,r')} \abs{f-a}^2\,d\mu \nonumber \\
&\ge \min_{a \in \bR}\int_{B_\Psi(x,r)} \abs{f-a}^2\,d\mu = \int_{B_\Psi(x,r)} \abs{f-f_{B_\Psi(x,r)}}^2\,d\mu. 
\end{align}

The Poincar\'e inequality  \hyperlink{pib}{$\operatorname{PI}(\beta)$} for $(\sX,d_\Psi,\mu,\sE,\sF^\mu)$ follows from \eqref{e:pcq2} and \eqref{e:pcq3}.
The converse is similar.\\
(b) Let $(\sX,d,\mu,\sE,\sF^\mu)$ satisfy \hyperlink{cs}{$\operatorname{CS}(\Psi)$} with constants $C_1,C_2,A \ge 1$.
Let $x \in X, r>0$ be arbitrary and let $r'>0, A'>1$ be as given in \eqref{e:pcq1}.
By \hyperlink{cs}{$\operatorname{CS}(\Psi)$}, there exists a cutoff function $\vp$ for $B(x,r') \subset B(x,Ar')$ such  that for any $u \in \sF^\mu \cap L^\infty$,
		\be \label{e:pcq4}
		\int_{B(x,Ar') \setminus B(x,r')} u^2 d\Gam(\vp,\vp) \le C_1 \int_{B(x,Ar') \setminus B(x,r')} d\Gam(u,u)
		+ \frac{C_2}{\Psi(x,r')} \int_{B(x,Ar') \setminus B(x,r')} u^2 \,d\mu.
		\ee 
Clearly by \eqref{e:pcq1}, $\vp$ is also a cutoff function for $B_\Psi(x,r) \subset B_\Psi(A' r)$. Since $\operatorname{supp} \Gam(\vp,\vp) \subset B(x,Ar') \setminus B(x,r')$, by \eqref{e:pcq1}, we have 
\be \label{e:pcq5}
\int_{B_\Psi(x,A'r) \setminus B_\Psi(x,r')} u^2 d\Gam(\vp,\vp) = \int_{B(x,Ar') \setminus B(x,r')} u^2 d\Gam(\vp,\vp).
\ee
Combining \eqref{e:pcq4}, \eqref{e:pcq5}, \eqref{e:pcq1}, and Proposition \ref{p:metric}(c), we obtain the cutoff energy inequality \hyperlink{csb}{$\operatorname{CS}(\beta)$} for $(\sX,d_\Psi,\mu,\sE,\sF^\mu)$. The converse is again similar.
\qed 

We will extend \hyperlink{cs}{$\operatorname{CS}(\Psi)$} to an inequality
for cutoff functions for $B(x,R) \subset B(x,R+r)$.
We will use the following elementary inequality involving energy measures.

\begin{lemma} \label{L:energy}
	Let $(\sE,\sF^\mu)$ be a regular Dirichlet form on $L^2(\sX,\mu)$ with energy measure $\Gamma(\cdot,\cdot)$. Then for any quasi-continous functions $f,\vp_1,\vp_2 \in \sF^\mu \cap L^\infty$, we have
	\[
	\int_{\sX} f^2 \, d\Gamma(\vp_1 \vee \vp_2,\vp_1 \vee \vp_2)
	 \le 	\int_{\sX} f^2 \, d\Gamma(\vp_1,\vp_1) + 	\int_{\sX} f^2 \, d\Gamma(\vp_2, \vp_2).
	\]
\end{lemma}
\proof 
Let $\vp_0= \vp_1 \vee \vp_2$. By \cite[Theorem 1.4.2(i),(ii)]{FOT}, we have 
$\vp_0 \in \sF^\mu, f^2 \in \sF^\mu$.
By \cite[last equation in p.206]{FOT} we have for each $j$
	\[
	\int_{\sX} f^2 \, d\Gamma(\vp_j,\vp_j)=  \lim_{t \downarrow 0} \frac{1}{t} \bE_{f^2.\mu} \left( (\vp_j(Y_t)- \vp_j(Y_0))^2 \right).
	\]
Here 
$\bE_{f^2.\mu}$ denotes the expectation where $Y_0$ has the distribution $f^2 \,d\mu$.
	Combining this with the elementary estimate,
	\[
	(\vp_0(Y_t)- \vp_0(Y_0))^2 \le \max_{i=1,2} 	(\vp_i(Y_t)- \vp_i(Y_0))^2 \le \sum_{i=1}^2(\vp_i(Y_t)- \vp_i(Y_0))^2,
	\]
	we obtain the desired inequality. 
\qed \\
The cutoff energy inequality  \hyperlink{cs}{$\operatorname{CS}(\Psi)$} has the following self improving property.

\begin{proposition} \label{P:selfimprove} 
(Cutoff energy inequality for all annuli) 
Let $(\sX,d,\mu,\sE,\sF^\mu)$ satisfy (MD) 
and \hyperlink{cs}{$\operatorname{CS}(\Psi)$} for some regular scale function  $\Psi$.
There exist $C_E$, $\gamma>0$ such that the following holds.
For all $R>0$, $r>0$, $x_0 \in \sX$ with $B_1=B(x_0,R)$, $B_2=B(x_0,R+r)$
and $U= B_2\setminus B_1$, there exists a cutoff function $\vp$
for $B_1 \subset B_2$ such that for any $f \in \sF^\mu \cap L^\infty$,
\be \label{e:modCS}
 \int_U f^2 d\Gam(\vp,\vp) \le \fract18 \int_U d\Gam(f,f)
  + C_E \left(\frac{R+r}{r}\right)^\gamma \frac{1}{\Psi(x_0,r)} \int_U f^2 d\mu.
\ee
\end{proposition}

\proof
Let $f \in \sF^\mu \cap L^\infty$.  
Let $A>1$, $C_1,C_2$ be the constants in \hyperlink{cs}{$\operatorname{CS}(\Psi)$}. 
Replacing $A$ by $\lceil A \rceil$ if necessary, we can assume that $A \in \bN$.
Let $n \ge 8(A+8)$ and cover $B(x_0,R+r)$ by balls $B_i= B(z_i, r/n), i \in I$ such that
$z_i\in B(x_0, R+r)$ and 
the balls $B(z_i, r/2n)$ are disjoint. Then using (MD) 
there exists a constant $N$ (which does not depend on $n$) such that any $y \in U$ 
is in at most $N$
of the balls $B_i^* = B(z_i, Ar/n)$. Let $U_i= B_i^* \setminus B_i$.

By \hyperlink{cs}{$\operatorname{CS}(\Psi)$}, there exists a cutoff function $\vp_i$ for $B_i \subset B_i^*$ such that
\begin{equation}\label{e:cut1}
\int_{U_i} f^2 \, d\Gamma(\vp_i,\vp_i) 
\le C_1 \int_{U_i} d\Gamma(f,f) + \frac{C_2}{\Psi(z_i,r/n)} \int_{U_i} f^2\, d\mu.
\end{equation}

Now let $2\le j \le n- A -1, j \in \bN$, and let
$I_j = \{ i \in I: z_i \in B(x_0,R+ jr/n)\}$. 
Set
$$ \psi_j = \max_{i \in I_j} \vp_i. $$
Then $\psi_j\equiv 1$ on $ B(x_0, R+ (j-1)r/n)$, and is zero outside
$B(x_0, R + (j+A) r/n)$. 
Thus $\psi_j$ is a cutoff function for $B(x_0, R+ (j-2)r/n) \subset B(x_0, R + (j+A+1) r/n)$.
We have $d(z_i,x_0) \le R+r$ for all $i \in I$, so using \eqref{e:Psireg1}
\be \label{e:cut1a}  \frac{\Psi(x_0, r)}{\Psi(z_i, r/n)}  
\le C \left(\frac{R+r}{r} \right)^{\beta_2-\beta_1} n^{\beta_2} . \ee

Let $V_j=  B(x_0, R + (j+A+1) r/n) \setminus B(x_0, R+ (j-2)r/n)$, so that $\supp\left(\Gamma(\psi_j,\psi_j)\right) \subset V_j$.

Let $h_j$ be a cutoff function for $\supp\left(\Gamma(\psi_j,\psi_j)\right) \subset V_j$.
By Lemma \ref{L:energy}
\begin{align}
\nn
  \int_{\sX} f^2 d\Gam(\psi_j, \psi_j) &=  \int_{\sX} f^2 h_j d\Gam(\psi_j, \psi_j) \\
   \label{e:cut1b}
  &\le \sum_{i \in I_j}  \int_{\sX} f ^2h_j d\Gam( \vp_i, \vp_i) \le \sum_{i \in I_j}  \int_{V_j} f ^2 d\Gam( \vp_i, \vp_i).
\end{align}

Now let
\bes \label{e:cut2}
 \vp = \frac{1}{ n-2A-4} \sum_{j=A+3}^{ n-A-2 }  \psi_j. 
\ees
Then $\vp$ is a cutoff function for $B(x_0,R) \subset B(x_0, R+r)$.
Since every point in $B(x_0,R+r)$ is in the support of at most $A+4$ of the energy measures $\Gam(\psi_j,\psi_j)$, by Cauchy-Schwarz inequality we have
\be \label{e:cut3}
 \int_{\sX} f^2 d\Gam(\vp,\vp) \le (A+4) (n-2A-4)^{-2} \sum_{j=A+3}^{n-A-2} \int_{\sX} f^2 d\Gam(  \psi_j, \psi_j).
\ee
 Combining \eqref{e:cut1b} and \eqref{e:cut3}, 
\bes 
\int_{\sX} f^2\, d\Gam(\vp,\vp) \le (A+4) (n-2A-4)^{-2} \sum_{j=A+3}^{n-A-2} \sum_{i \in I} \int_{V_j} f ^2 \,d\Gam( \vp_i, \vp_i).
\ees
Set $\wt I = \{ i \in I: \supp (\Gam(\vp_i,\vp_i)) \subset B(x_0,R+r) \setminus B(x_0,R) \}$.
If $A+2 \le j \le n-A-1$ and  $\supp\left( \Gam( \vp_i, \vp_i) \right) \cap V_j \neq \emptyset$ 
by the triangle inequality $\supp\left( \Gam( \vp_i, \vp_i) \right)  \subset B(x_0,R+r) \setminus B(x_0,R)$. 
Therefore it suffices to consider only the indices $i \in \wt I$ in \eqref{e:cut1b}.
Since for each $i$, $\supp(\Gam(\vp_i,\vp_i))$ intersects at most $4(A+4)$ different $V_j$'s we have, 
\be \label{e:cut4}
\int_{\sX} f^2\, d\Gam(\vp,\vp) \le 4(A+4)^2 (n-2A-4)^{-2}  \sum_{i \in \wt I} \int_{B(x_0,R+r) \setminus B(x_0,R)} f ^2 \,d\Gam( \vp_i, \vp_i).
\ee
Combining \eqref{e:cut4}, \eqref{e:cut1}, and \eqref{e:cut1a}, and using that every point 
is in at most $N$ different $B_i^*$, we obtain
\begin{align*}
 \int_{\sX} & f^2 d\Gam(\vp,\vp) 
   \\
  &\le \frac{4(A+4)^2}{(n-2A-4)^{2}} \left( C_1 \sum_{i \in \wt I } \int_{U_i}  d\Gam(f,f) 
      +  \frac{C_2 C n^{\beta_2} }{\Psi(x_0, r)}  
       \left(\frac{R+r}{r} \right)^{\beta_2-\beta_1}  \sum_{i \in \wt I}  \int_{U_i} f^2\, d\mu \right)\\
     &\le \frac{4N(A+4)^2}{(n-2A-4)^{2}} \left( C_1  \int_{U}  d\Gam(f,f) 
      +  \frac{C_2 C n^{\beta_2} }{\Psi(x_0, r)}  
       \left(\frac{R+r}{r} \right)^{\beta_2-\beta_1}  \int_{U} f^2\, d\mu \right).    
\end{align*}
Finally, we
choose $n$ large enough so that $4N(A+4)^2 (n-2A-4)^{-2} C_1 \le 1/8$. \qed

\begin{remark} \label{r:cs}
{\rm Note that this quite general argument enables us to deduce 
a cutoff energy inequality on arbitrary annuli from \hyperlink{cs}{$\operatorname{CS}(\Psi)$}. See \cite[Lemma 2.1]{MS2}. Further, if $\Psi(x,r)=\Psi(y,r)$ for all $x,y \in X$ and $r>0$, we can modify the proof by using \eqref{e:Psireg1} with $x=y$ instead of using \eqref{e:cut1a}, so that $\gamma=0$ in \eqref{e:modCS}.}
\end{remark}

\begin{definition} \label{d:cap-psi}
		Let $(\sX,d,\mu,\sE,\sF^\mu)$ be a MMD space and let $\Psi$ be a regular scale function.
		We say that $(\sX,d,\mu,\sE,\sF^\mu)$ satisfies the capacity estimate \hypertarget{cap}{$(\operatorname{cap})_\Psi$}, if there exists $\kappa \in (0,1)$ and $C>1$ such that for any ball $x \in \sX, r >0$, 
		\be   \tag*{$(\operatorname{cap})_\Psi$}
		C^{-1} \frac{\mu\left(B(x,r) \right)}{\Psi(x,r)} \le \Cpc_{B(x,r)}(B(x, \kappa r)) \le C \frac{\mu\left(B(x,r) \right)}{\Psi(x,r)}.
		\ee
			If $\Psi(x,r)=r^\beta$ for all $x \in \sX,r>0$, we denote the property 	 {$(\operatorname{cap})_\Psi$} by  \hypertarget{capb}{$(\operatorname{cap})_\beta$}. 
\end{definition}

We will now apply these results in the context of a change of measure on an MMD space.
Let  $(\sX,d,m,\sE,\sF^m)$ be a MMD space 
which satisfies the EHI and Assumptions \ref{A:green} and \ref{A:BG}.
Let $(\sE,\sF_e)$ denote the corresponding extended Dirichlet space (cf. \cite[Lemma 1.5.4]{FOT}), and 
$\mu$ be the measure constructed in Theorem \ref{T:meas}. 
By construction  $\mu$ is a positive Radon measure charging no set of capacity zero and possessing full support. 
Let $(\sE^\mu,\sF^\mu)$ denote the time changed Dirichlet space with respect to $\mu$. 
We have that $\sF^m=\sF_e \cap L^2(\sX,m)$,
$\sF^\mu=\sF_e \cap L^2(\sX,\mu)$
and $\sE^\mu(f,f)= \sE(f,f)$ for all $f \in \sF^\mu$ (Cf. \cite[p. 275]{FOT}). 
Moreover, the domain of the extended Dirichlet space is 
the same for both the Dirichlet forms  $(\sE,\sF^m,L^2(\sX,m))$ and $(\sE,\sF^\mu,L^2(\sX,\mu))$.

\begin{thm} \label{T:PI-CS}
$(\sX,d,m,\sE,\sF^m)$ be a length MMD space 
which satisfies the EHI and Assumptions \ref{A:green} and \ref{A:BG}.
Let $\mu$ be the measure constructed in Theorem \ref{T:meas}. 
Then the function $\Psi$ defined by $\Psi(x,0)=0$ and 
\be
\label{e:Psidef} \Psi(x,r) = \frac{\mu(B(x,r))}{ \Cpc_{B(x, r)}(B(x,r/8))}, \q r>0,
\ee
is a regular scale function.  Furthermore, the MMD space  $(\sX,d,\mu,\sE,\sF^\mu)$ satisfies  the Poincar\'e inequality \hyperlink{pi}{$\operatorname{PI}(\Psi)$}  and the cutoff 
energy inequality  \hyperlink{cs}{$\operatorname{CS}(\Psi)$}.
\end{thm}
\proof
By volume doubling and Corollary \ref{C:cap12}, there exists $C_2>0$ such that for all $r>0$ and for all $x,y \in \sX$ with $d(x,y) \le r$, we have
\be \label{e:Psi_comp}
C_2^{-1} \Psi(x,r) \le \Psi(y,r) \le C_2  \Psi(x,r).
\ee
If $R\le r$ the inequalities are immediate from property (b) in Theorem \ref{T:meas} and \eqref{e:Psi_comp}.
If $s<r<R$, then writing
\bes
\frac{ \Psi(x,r)}{\Psi(y,s)} =  \frac{ \Psi(x,r)}{\Psi(x,R)}  . \frac{ \Psi(y,R)}{\Psi(y,s)} .  \frac{ \Psi(x,R)}{\Psi(y,R)}, 
\ees
and bounding each of the three terms on the right using Theorem \ref{T:meas} and 
\eqref{e:Psi_comp} gives \eqref{e:Psireg1}. Thus $\Psi$ is a regular scale function.

Let $d_\Psi$ and $\beta>0$ be as given by Proposition \ref{p:metric}. 
Write $B_\Psi(\cdot,\cdot)$ for balls in the $d_\Psi$ metric.
We now show that  $(\sX,d_\Psi,\mu,\sE,\sF^\mu)$ satisfies \hyperlink{cap}{$(\operatorname{cap})_\beta$}.
By Lemma \ref{l:cann}, there exists $A>8, \kappa \in (0,1)$ such that for all $x \in X$, $r>0$, 
\[
B(x,s_1) \subset B_\Psi(x,\kappa r) \subset B(x,s_2) \subset B(x,8s_2) \subset B_\Psi(x, r) \subset B(x,As_1),
\]
for some $s_1,s_2>0$. 
By domain monotonicity of capacity, we have
\be \label{e:aic1}
\Cpc_{B(x,As_1)} \left( B(x,s_1) \right) \le \Cpc_{ B_\Psi(x, r) } \left( B_\Psi(x,\kappa r)  \right) \le \Cpc_{B(x,8s_2)} \left( B(x,s_2) \right). 
\ee
By Proposition \ref{p:metric}(c) and the regularity of $\Psi$, $s_1$ and $s_2$ are both 
comparable with $\Psi(x,s_1) \asymp \Psi(x,s_2) \asymp r^\beta$. Therefore by (VD), Lemmas \ref{L:cap1}, \ref{L:domcap}, \eqref{e:ggrowth}, and \eqref{e:aic1}, we have  
\be  \label{e:aic2}
  \Cpc_{B(x,As_1)} \left( B(x,s_1) \right) \asymp  \Cpc_{B(x,8s_2)} \left( B(x,s_2) \right) \asymp \frac{\mu(B(x,s_2))}{\Psi(x,s_2)} \asymp \frac{\mu(B_\Psi(x,r))}{r^\beta}.
\ee
Combining \eqref{e:aic1} and \eqref{e:aic2}, we have that $(\sX,d_\Psi,\mu,\sE,\sF^\mu)$ satisfies \hyperlink{cap}{$(\operatorname{cap})_\beta$}.
By Lemma \ref{l:ehiqs} and Proposition \ref{p:metric}(b), $(\sX,d_\Psi,\mu,\sE,\sF^\mu)$ satisfies the EHI. 

By Remark \ref{R:up} the space $(\sX,d_\Psi)$ is uniformly perfect, and the measure $\mu$
on $(\sX,d_\Psi)$ satisfies (RVD). 
Therefore by \cite[Theorem 1.2]{GHL}, since $(\sX,d_\Psi,\mu,\sE,\sF^\mu)$  satisfies  
the EHI and \hyperlink{cap}{$(\operatorname{cap})_\beta$}, it satisfies 
\hyperlink{pib}{$\operatorname{PI}(\beta)$}  and \hyperlink{csb}{$\operatorname{CS}(\beta)$}. 
We now conclude using Lemma \ref{l:picsqs}. \qed

\begin{thm} \label{T:main-new}
Let $(\sX,d)$ be a complete, locally compact, length metric space with a strongly local 
regular Dirichlet form $(\sE,\sF^m)$ on $L^2(\sX,m)$ which satisfies
Assumptions \ref{A:green} and \ref{A:BG}. 
The following are equivalent: \\
(a)  $(\sX,d,m,\sE,\sF^m)$ satisfies the EHI. \\
(b) There exists a doubling Radon measure $\mu$ on $(\sX,d)$ 
which is mutually absolutely continuous with respect to $m$, and a regular scale function $\Psi$,
 such that the time-changed MMD space  $(\sX,d,\mu,\sE,\sF^\mu)$ satisfies the Poincar\'e inequality \hyperlink{pi}{$\operatorname{PI}(\Psi)$}  and the cutoff 
energy inequality  \hyperlink{cs}{$\operatorname{CS}(\Psi)$}.\\
(c) There exists a doubling Radon measure $\mu$ on $(\sX,d)$ 
which is mutually absolutely continuous with respect to $m$, a  metric $d_\Psi$ on $\sX$ that is quasisymmetric to $d$, and $\beta>0$, such that the time-changed MMD space   $(\sX,d_\Psi,\mu,\sE,\sF^\mu)$ satisfies Poincar\'e inequality \hyperlink{pib}{$\operatorname{PI}(\beta)$}  and the cutoff 
energy inequality  \hyperlink{csb}{$\operatorname{CS}(\beta)$} for some $\beta >0$.
\end{thm}
\proof 
(a) $\Rightarrow $(b) This follows from Theorem \ref{T:PI-CS}.
\\
(b)$\Rightarrow$(c) Let $d_\Psi$ and $\beta>0$ be as given by Proposition \ref{p:metric}. 
Quasisymmetry of $d_\Psi$ follows from Proposition \ref{p:metric}(b). 
Then \hyperlink{pib}{$\operatorname{PI}(\beta)$} and \hyperlink{csb}{$\operatorname{CS}(\beta)$}  for $(\sX,d_\Psi,\mu,\sE,\sF^\mu)$  follow from Lemma \ref{l:picsqs}. \\
(c)$\Rightarrow$(a)
By Remark \ref{R:up} $(\sX,d)$ are therefore $(\sX,d_\Psi)$ are uniformly perfect. Thus
$\mu$ satisfies (RVD). 
 By Proposition \ref{P:selfimprove} and Remark \ref{r:cs}, 
we obtain the condition (CSA) in \cite{GHL}. Then by  \cite[Theorem 1.2]{GHL}, we obtain EHI for $(\sX,d_\Psi,\mu,\sE,\sF^\mu)$. Since $d_\Psi$ and $d$ are quasisymmetric, the desired EHI follows from Lemma \ref{l:ehiqs}.
\qed

{\smallskip\noindent {\em Proof of Theorem \ref{T:eeprime}.}}
The relation $\sE \asymp \sE'$ implies that the energy 
measure $d\Gam'(f,f)$ for $\sE'$ satisfies
\be \label{e:stabEN}
 C^{-1} d\Gam(f,f) \le d\Gam'(f,f) \le C d\Gam(f,f) \hbox{ for all } f \in \sF
\ee
-- see \cite[Proposition 1.5.5(b)]{LJ}. 
This implies stability of Poincar\'e and cutoff energy inequalities under such perturbations.
Therefore, 
the desired EHI follows from stability of property (b) in Theorem \ref{T:main-new} (or alternatively (c)).
\qed\\

\sms
We remark that the cutoff energy inequality in Theorems 
\ref{T:main-new} and \ref{T:PI-CS} could be 
replaced by the slightly weaker generalized capacity estimate given in \cite{GHL}.
\begin{remark} \label{r:qs} {\rm
(1) The approach using quasisymmetry  given in this section implicitly contains an alternate proof to the main results in \cite{Bas}. \\
(2) Theorem \ref{T:main-new} shows that,
after suitable transformations of measure and metric, the 
 stability of EHI follows from the stability of the $\operatorname{PHI}(\beta)$ -- see \cite[p. 485 and Definition 2.1(d)]{BBK} for the definition of $\operatorname{PHI}(\beta)$. 
It is known that the index $\beta \ge 2$  -- see \cite[p. 252]{Hin}.
One might ask if we can improve Theorem \ref{T:main-new}(c)  to  obtain $\operatorname{PHI}(2)$.
A paper in preparation \cite{KM} shows that this is not possible in general, but
on the other hand 
the Sierpinski gasket provides a non-trivial  example where this is possible -- see \cite{Ki}. 
See \cite[Section 9]{Ka} for further discussion on this problem.\\
(3) The constant $\beta >0$ in Theorem \ref{T:main-new} can be made arbitrarily 
large by a `snowflake transform' of the metric 
 $d_\Psi \mapsto d_\Psi^\varepsilon$, where $\varepsilon \in (0,1)$. 
 We can ask how small $\beta$ can be.
Recall that a \emph{conformal gauge} on a set $\sX$ is a maximal collection of metrics on $\sX$ such that each pair of metrics from the collection are quasisymmetric. By analogy with conformal Hausdorff dimension (see \cite[Defintion 2.2.1]{MT} or \cite[pg. 121]{Hei}), we can define the  \emph{conformal walk dimension} 
of a MMD space $(\sX,d,m,\sE,\sF^m)$ as the infimum of all $\beta$ such that there exists a 
quasisymmetric metric $d_\Psi$ and a Revuz measure $\mu$ with full support such that the time 
changed MMD space $(\sX,d_\Psi,\mu,\sE,\sF^\mu)$ satisfies $\operatorname{PHI}(\beta)$. 
The conformal walk dimension is always at least 2, and 
by Theorem \ref{T:main-new} it is finite if and only if the space satisfies EHI. 
This raises the following questions:  Can the conformal walk dimension be finite and strictly greater than 2? Is the infimum in the definition of  conformal walk dimension always attained? \\
(4) By \cite[Theorem 1.2]{GHL} the modified space $(\sX,d_\Psi,\mu,\sE,\sF^\mu)$ satisfies
heat kernel upper and lower bounds -- see \cite{GHL} for details. \\
(5) The classical parabolic Harnack inequality $\operatorname{PHI}(2)$ implies that vector space of harmonic functions with fixed polynomial growth is finite dimensional  \cite[Theorem 0.7]{CM}. 
This result of Colding and Minicozzi settled a conjecture of Yau on manifolds with non-negative Ricci curvature.
This result was extended by P. Li \cite[Theorem 1]{Li} to spaces satisfying a mean value inequality for harmonic functions with respect to a doubling measure. 
This theorem of Li along with our doubling measure $\mu$ in Theorem \ref{T:meas} implies that the vector space of harmonic functions with fixed polynomial growth is finite dimensional on any space satisfying the EHI. 
Note that one cannot directly use \cite[Theorem 1]{Li} to obtain the above result because there are manifolds that satisfy EHI but whose Riemannain measure is not doubling. 
}\end{remark}

\section{Examples: Weighted Riemannian manifolds and graphs}
\label{sec:ex}

In this section we return to our two main examples,  
and give sufficient conditions 
for these spaces to satisfy the local regularity hypotheses
\ref{A:green} and \ref{A:BG}.

We first recall some standard definitions in Riemannian geometry. 
Let $(\sX,g)$ be a Riemannian manifold,
and $\nu$ and  $\nabla$  denote the Riemannian measure  and the Riemannian gradient respectively. In local coordinates $(x_1,x_2,\ldots,x_n)$, we have
	\[
	\nabla f = \sum_{i,j=1}^n g^{i,j} \frac{\partial f}{\partial x_i} \frac{\partial}{\partial x_j}, \hspace{4mm} d\nu = \sqrt{\det g(x)}\,dx,
	\] 
	where $\det g$ denotes the determinant of the metric tensor $(g_{i,j})$ and $(g^{i,j})=(g_{i,j})^{-1}$ is the co-metric tensor. For a function $f \in \sC^\infty(\sX)$, we denote the length of the gradient by $\abs{\nabla f}= \left( g(\nabla f,\nabla f)\right)^{1/2}$.
The Laplace-Beltrami operator $\Delta$ is given in local coordinates by
	\[
	\Delta  = \frac{1}{ \sqrt{\det g}} \sum_{i,j} \frac{\partial}{\partial x_i} \left( g^{i,j} \sqrt{\det g} \frac{\partial}{\partial x_j} \right).
	\]

A {\em weighted manifold} $(\sX,g,\mu)$ is a Riemmanian manifold 
$(\sX,g)$ endowed with a measure $\mu$ that has a smooth (strictly) 
positive density $w$ with respect to $\nu$. 
Let $w$ be the smooth function such that \[d \mu= w d \nu.\] 
On the weighted manifold $(M,g,\mu)$, one associates a {\em weighted Laplace operator} $\Delta_\mu$ given  by
\[
\Delta_\mu f= \Delta f + g\left( \nabla \left( \ln w \right), \nabla f  \right), \mbox{ for all $f \in \sC^\infty(\sX)$}.
\]
We say that the weighted manifold $(M,g, \mu)$ has {\em controlled weights} if the function $w$ defined above satisfies 
\[
\sup_{x,y \in \sX : d(x,y) \le 1} \frac{w(x)}{w(y)} < \infty,
\]
where $d$ is the Riemannian distance function.
The corresponding Dirichlet form on $L^2(\sX,\mu)$ is given by
\[
\sE(f_1,f_2)= \int_{\sX} g(\nabla f_1, \nabla f_2)\,d\mu, \q f_1,f_2 \in \sF,
\] 
where $\sF$ is the weighted Sobolev space of functions in $L^2(\sX,\mu)$ whose distributional gradient is also in $L^2(\sX,\mu)$.  We refer the reader to Grigor'yan's survey \cite{Gri06} for details of the construction of the heat kernel, Markov semigroup and Brownian motion on weighted manifolds for motivation, as well
as applications. 

Our second example is weighted graphs. 
Let $\bG = (\bV,E)$ be an infinite graph, such that each vertex 
$x$ has finite degree. For $x \in V$ we
write $x \sim y$ if $\{x,y\} \in E$. For $D \subset \bV$ define
$$ \pd D =\{ y \in D^\compl: y \sim x \hbox{ for some } x \in D \}. $$
We define a metric on $V$ by taking $d(x,y)$ to be the
length of the shortest path connecting $x$ and $y$.
We define balls by
$$ B_d(x,r) = \{ y\in \bV: d(x,y) < r\}. $$

Let $w: E \to (0,\infty)$ be a function
which assigns weight $w_e$ to the edge $e$. We write
$w_{xy}$ for $w_{\{x,y\}}$, and define
\be \label{e:wx}
w_x = \sum_{y \sim x} w_{xy}. 
\ee
We extend $w$ to a measure on $\bV$ by setting $w(A) =\sum_{x \in A} w_x$.
We call $(\bV,E,w)$ a {\em weighted graph}. An {\em unweighted graph} has
$w_e \equiv 1$.

The Dirichlet form associated with this weighted graph is given by taking
$$ \sE_\bG(f,f) = \half \sum_{x} \sum_{y \sim x} w_{xy} (f(y)-f(x))^2, $$
with domain $\sF =\Sett{ f \in L^2(\bV,w)}{ \sE_\bG(f,f) < \infty}$.
We define the Laplacian on $\bG$ by setting
$$ \Delta_\bG f(x) =  \frac{1}{w_x} \sum_{y \sim x}  w_{xy} (f(y)-f(x)). $$
We say that a function $h$ is {\em harmonic} on a set $D \subset \bV$ if
$\Delta_\bG h(x) =0$ for all $x \in D$. 
(Note that for $\Delta_\bG h(x)$ to be defined for $x \in D$ we need $h$
to be defined on the set $D \cup \pd D$.)

The statement of the elliptic Harnack inequality for a weighted
graph is analogous to the EHI for a MMD space. 
We say $\bG=(V,E,w)$ satisfies the EHI is there exists $C_H<\infty$ such that
if $x_0 \in \bV$, $R\ge 1$, and $h: B(x_0, 2R+1)\to \bR_+$ is harmonic in
$B(x_0, 2R)$ then
$$ \sup_{B_d(x_0,R)} h \le C_H   \inf_{B_d(x_0,R)} h. $$

The {\em cable system of a weighted graph}
gives a natural embedding of a graph in a connected metric length space.
Choose a direction for each edge $e \in E$, let $(I_e, e \in E)$ be a collection
of copies of the open unit interval, and set
$$ \sX = \bV \cup (\cup_{e \in E} I_e) . $$
(Following \cite{V} we call the sets $I_e$ {\em cables}).
We define a metric $d_c$ on $\sX$ by using Euclidean distance on each cable.
If $x \in \bV$ and $e=(x,y)$ is an oriented edge, we set
$d_c(x,t) = 1-d_c(y,t) = t$ for $t \in I_e$. We then extend $d_c$ to a metric
on $\sX$; note that this agrees with the graph metric for $x,y \in \bV$.
We take $m$ to be the measure on $\sX$ which assigns zero mass to 
points in $\bV$, and mass $w_e |s-t|$ to any interval $(s,t) \subset I_e$.
For more details on this construction see \cite{V,BB3}.

We say that a function $f$ on $\sX$ is piecewise differentiable
if it is continuous at each vertex $x \in \bV$, is differentiable on each cable, and
has one sided derivatives at the endpoints. Let $\sF_{0}$ be the set
of piecewise differentiable functions $f$ with compact support.
Given two such functions we set
$$ d\Gam(f,g)(t) =  f'(t)g'(t) m(dt). $$
(While the sign of $f'$ and $g'$ depends on the orientation of the
cable this does not affect their product.) 
We then define
\begin{align*}
 \sE(f,g) = \int_{\sX}  d\Gam(f,g)(t), \q f,g \in \sF_0,
\end{align*}
and take $\sF$ to be the completion of $\sF_0$ with respect to the norm
$$ ||f||_{\sE_1} = \Big ( \int f^2 dm + \sE(f,f) \Big)^{1/2}. $$
We extend $\sE$ to $\sF$, and it is straightforward to verify 
that $(\sE,\sF)$ is a closed regular strongly local Dirichlet form. 
We call $(\sX,d_c,m,\sE,\sF)$ the {\em cable system} of the graph $\bG$.
We define harmonic functions for the cable system as in Section \ref{sec:intro}.

We remark that (up to a constant time change) 
the associated Hunt process $X$ behaves like a Brownian
motion on each cable, and like a `Walsh Brownian motion' (see \cite{W})
at each vertex: starting at $x$ it makes excursions along the cable
$I_{\{x,y\}}$ at rate proportional to $w_{xy}/w_x$. 

There is a natural bijection between harmonic functions 
on the graph $\bG$ and the cable system $\sX$. 
If $h$ is harmonic on a domain $D \subset \sX$
then $h|_{\bV}$ satisfies $ \Delta_\bG h(x)=0$ for any $x \in \bV$
such that $B(x,1) \subset D$. Conversely  let
$D_0 \subset \bV$, and suppose that $h : D_0 \cup \pd D_0 \to \bR$
is $\bG$-harmonic. Let $D$ be the open subset of $\sX$ which consists
of $D_0$ and all cables with an endpoint in $D_0$. Define $\ol h$
by setting $\ol h(x) = h(x)$, $x \in D_0 \cup \pd D_0$, and taking $\ol h$
to be linear on each cable. Then $\ol h$ is harmonic on $D$.

\begin{definition}
{\rm We say that $\bG$ has {\em controlled weights} if there exists
$p_0>0$ such that
\be
  \frac{w_{xy}}{w_x} \ge p_0 \, \hbox{ for all } x\in \bV, \; y \sim x.
\ee
} \end{definition}

This is called the $p_0$ condition in \cite{GT}. Note that it implies
that vertices have degree at most $1/p_0$, so that an unweighted graph
satisfies controlled weights if and only if the vertex degrees are uniformly bounded.

\begin{lemma}
	Let $(\sX,d,\mu,\sE,\sF)$ be the cable system of a weighted graph 
	$\bG=(\bV, E, w)$. 
	If $\sX$ satisfies the EHI with constant $C_H$ then
	$\bG$ has  controlled weights.
\end{lemma}

\proof (By looking at a linear (harmonic) function in a single cable 
we have that $C_H \ge 3$.) Let $x_0 \in \bV$ and let $x_i$, $i=1, \dots n$
be the neighbours of $x_0$. Let $r < \half$, and $y_i, z_i$ be the points
on the cable $\gam(x_0, x_i)$ with $d(x_0,y_i)=r$,  $d(x_0,z_i)=2r$.
Set $p_j = w_{x_0,x_j}/w_x$. 

Let $D= B(x_0,2r)$ and $h_j$ be the harmonic function in $B(x_0,2r)$
with $h_j(z_i)=\delta_{ij}$.
We have $h_j(x_0) = p_j$, $h_j(y_i) = \half p_j$ if $i \neq j$ and
$h_j(y_j) = \half(1 + p_j)$. So using the EHI with $i\neq j$
$$  2 h(y_j) = 1 + p_j \le 2 C_H h(y_i) = C_H p_j, $$
which leads to the required lower bound on $p_j$.
\qed

\begin{remark}
	{\rm  See \cite{B1} for an example which shows that the EHI for a weighted
		graph, as opposed to its cable system, does not imply controlled weights.}
\end{remark}

It is straightforward to verify
\begin{lemma}
	Let $\bG$ have controlled weights. 
	The EHI holds for $\bG$ if and only if it holds for the associated cable system.
\end{lemma} 

We conclude this section by showing that a large class of weighted manifolds and cable systems satisfy our local regularity
hypotheses \BG. To this end, we  introduce a local parabolic Harnack inequality which turns out to be strong enough to imply \BG.

\begin{definition} 
We say a MMD space  $(\sX,d,\mu,\sE,\sF)$ satisfies the local parabolic Harnack inequality 
$(\operatorname{PHI(2)})_{\loc}$, if there exists $R>0,C_R >0$
	such that for all $x \in \sX$, $0< r \le R$, any non-negative weak solution 
	$u$ of $(\partial_t+ \sL)u=0$ on $(0,r^2) \times B(x,r)$ satisfies
	\begin{equation*} \label{phi2-loc}
	\sup_{(r^2/4,r^2/2) \times B(x,r/2)} u \le C_R \inf_{(3r^2/4,r^2) \times B(x,r/2)} u;
	\tag*{$(\operatorname{PHI(2)})_{\loc}$}
	\end{equation*} 
	here $\sL$ is the generator corresponding to the  Dirichlet form $(\sE,\sF, L^2(\sX,\mu))$.
\end{definition}

\begin{lemma}  \label{l:PHIloc} 
(a) Let $(\sM,g,w)$ be a weighted Riemannian manifold with controlled weights 
such that $(\sM,g)$ is quasi-isometric to a manifold with Ricci curvature 
bounded below. Then $(\sM,g,w)$ satisfies $(\operatorname{PHI(2)})_{\loc}$.\\
(b) Let $\bG=(\bV, E, w)$ be a weighted graph with controlled weights.
Then its cable system satisfies $(\operatorname{PHI(2)})_{\loc}$. 
\end{lemma}

\proof (a)
If $(\sM',g')$ has Ricci curvature bounded below then  $(\sM',g')$ satisfies
$(\operatorname{PHI(2)})_{\loc}$ by the Li-Yau estimates. 
By [HS, Theorem 2.7], the property $(\operatorname{PHI(2)})_{\loc}$ is stable under quasi isometries 
and under introducing controlled weights. \\
(b) By taking $R<1$ this reduces to looking at either a single cable (i.e. an interval)
or a finite union of cables. See \cite{BM} for more details. \qed

\begin{lemma}\label{l:phi-bg}
Let $(\sX,d,m,\sE,\sF)$ be a MMD space that satisfies 
\ref{phi2-loc}.  Then $(\sX,d,m,\sE,\sF)$ satisfies Assumption \ref{A:green} and \BG.
\end{lemma}

\proof
We refer the reader to \cite{BM} for the proof of Assumption \ref{A:green}.

By \cite[Theorem 2.7]{HS} the heat kernel on this space satisfies 
a two sided Gaussian bound at small time scales. 
These imply volume doubling property at small scales.

Using the heat kernel upper bounds given in \cite[Lemma 3.9]{HS}, 
we obtain the following Green's function upper bound. 
There exists $A >1$, $a \in (0,1), C_0, r_0 >0$ such that
for all $x \in \sX, r \in (0,r_0)$ and for all $y \in B(x,Ar)$ such that $d(x,y) = ar$, we have
\[
g_{B(x,Ar)}(x,y) \le C_0 \frac{r^2}{m(B(x,r))}.
\]
A matching lower bound 
follows from \cite[Lemmas 3.7 and 3.8]{HS},  after adjusting $r_0,a$ if necessary.

Clearly, \ref{phi2-loc} implies a local EHI for small scales. By using the local EHI along with the results in Section 2 (see Remark \ref{r-local}), there exists $r_0,C_1>0$ such that
\[
C_1^{-1} \frac{m(B(x,r))}{r^2}\le \Cpc_{B(x,8r)}(B(x,r)) \le \frac{m(B(x,r))}{r^2}, \q \forall x \in \sX, \forall r \in (0,r_0).
\]
This implies \eqref{e:eois} with $\gamma_2=2$. Hence \BG \hspace{1mm}follows.
\qed

{\smallskip\noindent {\em Proof of Theorem \ref{T:mainconsq}.}}
Assumption \ref{A:BG} follows from Lemma \ref{l:PHIloc} and \ref{l:phi-bg}. Assumption \ref{A:green} follows from \cite{BM}. The conclusions now follow from Theorem \ref{T:eeprime}.
\qed

\section{Stability under rough isometries} \label{sec:qi}

As well as stability of the EHI under bounded perturbation of weights,
our results also imply stability under rough isometries.

\begin{definition}
{\rm 
For each $i=1,2$, let $(\sY_i,d_i,\mu_i)$ be either 
a metric measure space or a weighted graph. 
A map
$\vp: \sY_1 \to \sY_2$ is a {\sl rough  isometry} 
if there exist constants $C_1>0$ and $C_2, C_3>1$ such that 
\begin{align}
\label{e:qi1}
X_2 &=\bigcup_{x\in X_1} B_{d_2}(\vp(x),C_1), \\
\label{e:qi2}
 C_2^{-1} (d_1(x,y)-c_1) &\le d_2(\vp(x),\vp(y))\le  C_2 (d_1(x,y)+c_1),  
\hbox{ for } x  \in \sY_1,  \\
\label{e:qi3}
C_3^{-1} \mu_{1}(B_{d_1}(x,C_1)) &\le  \mu_2(B_{d_2}(\vp(x),c_1)) 
 \le  C_3 \mu_{1}(B_{d_1}(x,C_1)) \hbox{ for } x,y  \in \sY_1.
 \end{align}
If there exists a rough isometry between two spaces they are
said to be {\sl roughly isometric}. (One can check this is an
equivalence relation.) 
}\end{definition}

This concept was introduced by Gromov \cite{Gro} (under the name quasi isometry) in the context of groups, 
and Kanai \cite{Kan} (under the name rough isometry) for metric spaces;
in both cases they just required the conditions \eqref{e:qi1} and \eqref{e:qi2}.
The condition \eqref{e:qi3} is a natural extension when one treats measure
spaces -- see \cite{CS}
 and \cite{BBK}.

If two spaces are roughly isometric then they have similar large scale structure. 
However, as the EHI  implies some local regularity, we
need to impose some local regularity on the spaces in the class we consider.

\begin{definition}
We say a MMD space satisfies a local EHI (denoted $\operatorname{EHI}_{\loc}$)
if there exists $r_0 \in(0,\infty)$ and $C_L<\infty$ such that whenever $2r<r_0$,
$x \in \sX$ and $h$ is a nonnegative harmonic function on $B(x,2r)$ then
$$ \esssup_{B(x,r)} h \le C_L  \essinf_{B(x,r)} h . $$
\end{definition}

\begin{remark}
{\rm 
An easy chaining argument shows that if $\sX$ satisfies EHI$_{\loc}$
with constants $r_0$ and $C_L$, then for any $r_1>r_0$ there exists
$C_L'=C_L(r_1)$ such that $\sX$ satisfies EHI$_{\loc}$
with constants $r_1$ and $C_L'$.
} \end{remark}

\begin{definition}
{\rm  Let $\sX=(\sX,d,m,\sE,\sF)$  be a MMD space.  
We say $\sX$ satisfies {\em local regularity} (LR) if there exists
$r_0 \in (0,1)$, $C_L<\infty$ such that the following conditions hold: \\
(B1) $\sX$ satisfies \BG.  \\
(B2)  The Green's function and operator satisfies Assumption \ref{A:green}. \\
(B3) $\sX$ satisfies $\operatorname{EHI}_{\loc}$ with constants $r_0$ and $C_L$.\\
(B4) There exists $C_0>0$ such that for all $x_0 \in \sX$ and for all $r \in (0,r_0)$, there exists a cut-off function $\vp$
for $B(x_0,r/2) \subset B(x_0,r)$ such that
\[
\int_{B(x_0,r)} d \Gamma(\vp,\vp) \le C_0 m(B(x_0,r)).
\]
}\end{definition}

The final  condition (B4) links $m$ with the energy measure 
$d\Gam(\cdot, \cdot)$ at small length scales.

\begin{lemma}  \label{e:L-LR}
(a) Let $(\sM,g,w)$ be a weighted Riemannian manifold with controlled weights such that $(\sM,g)$ is quasi-isometric to a manifold with Ricci curvature bounded below. Then $(\sM,g,w)$ satisfies (LR).\\
(b) Let $\bG=(\bV, E, w)$ be a weighted graph with controlled weights.
Then its cable system satisfies (LR).
\end{lemma}

\proof Properties (B1)--(B3) all follow from Lemma \ref{l:PHIloc} and \ref{l:phi-bg}.
For (B4) it is sufficient to look at a cutoff function $\vp(x)$ which is piecewise linear 
in $d(x,x_0)$.
\qed

Our main theorem concerning stability under rough isometries 
is the following.

\begin{thm}[Stability under rough isometries] \label{T:qistab}
Let $\sX_i=(\sX_i,d_i,m_i,\sE_i,\sF_i)$, $i=1,2$  be MMD spaces which satisfy
(LR). Suppose that $\sX_1$ satisfies the EHI, and $\sX_2$ is roughly isometric 
to $\sX_1$. Then $\sX_2$ satisfies the EHI.
\end{thm}
{\smallskip\noindent {\em Sketch of the proof.}}
The basic approach  goes back to the seminal works of Kanai \cite{Kan,Kan2,Kan3};
see \cite{CS,HK,BBK} for further developments.

We use the characterization of EHI
in  Theorem \ref{T:main-new}, and transfer functional inequalities 
and volume estimates from one space to the other. 
A  key step of this transfer is carried out by a discretization procedure using weighted graphs.

We can approximate an MMD space $(\sX,d,m,\sE,\sF^m)$ 
by a weighted graph as follows. 
For a small enough $\eps$, we choose an $\eps$-net $\bV$ 
of the MMD space $(\sX,d,m,\sE,\sF^m)$, that is a maximal $\eps$-separated subset of $X$.
The set $\bV$ forms the vertices of a graph whose edges 
$E$ are given by $u \sim v$ if and only if $d(u,v) \le 3 \eps$. 
Define weights by 
$w_{uv} = m(B(u,\eps))+ m(B(v,\eps))$ if $\set{u,v} \in E$.
(Many other choices are possible.)
We then define $w_x$ as in \eqref{e:wx} and hence obtain a measure $w$ on $\bV$.
It is easy to verify that the metric measure spaces $(\sX,d,m)$
 and $(\bV,E,w)$ are  roughly isometric.
 
The next step is to  transfer functions between MMD space and its net. 
This transfer of functions has the property that the norms and energy measures are 
 comparable on balls (up to constants and linear scaling of balls), 
 which in turn implies that functional inequalities such as the 
 Poincar\'e inequality and cutoff energy inequality can be transferred between a MMD space and 
 its net. Using the notation of \cite{Sal04}, we denote by  $\rst$ a ``restriction map" that takes a function $f:\sX \to \mathbb{R}$ on the MMD space  to a function 
 $\rst(f): \bV \to \mathbb{R}$ on the graph defined by 
 \[
\rst(f)(v) = \frac{1}{m(B(v,\eps))} \int_{B(v,\eps)} f(y)\,m(dy), \hbox{ for } v \in \bV.
 \]
 Similarly,  we denote by  $\ext$ an ``extension map" that takes a function 
 $f:\bV \to \mathbb{R}$ on the net  to a function $\ext(f): \sX \to \mathbb{R}$ on the MMD space defined by 
 \[
\ext(f)(x) = \sum_{v \in V} f(v) \chi_{v}(x), 
 \]
 where $(\chi_v)_{v \in V}$ is a `nice' 
 partition of unity on $\sX$ indexed by the vertices of the net $V$ satisfying the following properties:
 \begin{enumerate}[(i)]
 \item $\sum_{v \in V} \chi_v = 1$.
 \item There exists $c \in (0,1)$ such that  $\chi_v \ge c$ on $B(x,\eps/2)$ for all  $v \in V$.
 \item $\chi_v \equiv 0$ on $B(v,2 \eps)^\compl$ for all $v \in V$.
 \item There exists $C>0$ such that $\chi_v \in \sF^m$  and  $\sE(\chi_v,\chi_v) \le 
 C m(B(x,\eps))$ for all $v \in V$.
 \end{enumerate}
 The maps $\rst$ and $\ext$ are (roughly) inverses of each other,
 and they preserve norms and energy measures on balls. Therefore volume doubling, 
the Poincar\'{e} inequality,  and the cutoff energy inequality can be transferred 
between a MMD space and its net. 
 
 A difficulty that is not present in the previous settings 
 in \cite{CS,HK, BBK} arises from the change of measure in the
 characterization of the EHI in Theorem \ref{T:main-new}. 
 This change of symmetric measure does not affect the energy measures in the cutoff energy 
 and Poincar\'e inequalities. However the integrals on the left side of the Poincar\'e
 inequality, and the final integral in the cutoff energy inequality involve
 the measure measure $\mu$  constructed in Theorem \ref{T:meas}.
Let $g$ be such that  $d \mu =g dm$. 
The integrals
for the cutoff energy and Poincar\'{e} inequalities on
the net the  are taken with respect to the measure $\rst(g)\,d\mu$. 
It is easy to verify using \eqref{e:m03} that the 
metric measure spaces $(\sX,d,\mu)$ and the net equipped with the measure $\rst(g)\,dw $ are roughly 
isometric, and therefore integrals with respect to the measures 
$g \,dm$ and $\rst(g)\,dw$ are comparable on balls. 

Thus, starting with the space $\sX_1$ we take $g_1 = d\mu_1/dm_1$, where $\mu_1$
is the measure given by Theorem \ref{T:meas}. Write $\bV_i$ for the nets for $\sX_i$, $i=1,2$.
We take $\wt g_1 = \rst(g_1)$, and then transfer $\wt g_1$ to a function $\wt g_2$ on 
$\bV_2$ using the rough isometry between $\bV_1$ and $\bV_2$. The function
$g_2 = \ext(\wt g_2)$ then gives a measure $d\mu_2= g_2 dm_2$ on $\sX_2$.
As in \cite{CS,HK, BBK} we can then transfer the cutoff energy and Poincar\'{e} inequalities
across this chain of spaces, and deduce that the space 
$(\sX_2,d_2,\mu_2,\sE_2,\sF_2)$ satisfies the conditions in Theorem \ref{T:main-new}(b), 
and therefore satisfies the EHI.
\qed

{\smallskip\noindent {\em Proof of Theorem \ref{T:mainconsqG}.}}
This is a direct consequence of Lemma \ref{e:L-LR} and Theorem \ref{T:qistab}.
\qed

We conclude this paper by suggesting a characterization of the EHI in terms
of capacity, or equivalently effective conductance. 
Let $D$ be a bounded domain in $\sX$.
As in \cite{CF} we can define a reflected Dirichlet space $\wt \sF_D$; the associated
diffusion $\wt X$ is the process $X$ reflected on (a) boundary of $D$. 
(For the case of manifolds or graphs this reflected process can be constructed in a 
straightforward fashion). For disjoint subsets $A_1$, $A_2$ of $D$ define 
$$ \Ceff(A_1,A_2;D) = \inf\{ \sE_D(f,f): f|_{A_1}=1, f|_{A_2}=0, f \in \wt \sF_D \}. $$
Let $\sD(x_0,R)=\{ (x,y) \in B(x_0,R): x,y \in B(x_0,R/2), d(x,y) \ge R/3\}$.
As in \cite{B1} we say that $(\sX,d,m,\sE,\sF)$ satisfies the {\em dumbbell condition} 
if there exists $C_D$ such that for all $x_0\in \sX$, $R>0$ we have, writing $D=B(x_0,R)$,
\bes
 \sup_{(x,y)\in \sD(x_0,R) } \Ceff(B(x,R/8), B(y,R/8);D) 
  \le C_D  \inf_{(x,y)\in \sD(x_0,R) } \Ceff(B(x,R/8), B(y,R/8);D).
\ees
\cite{B1} asked if the dumbbell condition characterizes EHI.
However G. Kozma \cite{Ko} remarked that a class of spherically symmetric trees satisfy 
the dumbbell condition, but fail to satisfy EHI. These trees also fail to satisfy (MD). 
We can therefore modify the question in \cite{B1} as follows.

\begin{problem} {\rm
Let $(\sX,d,m,\sE,\sF^m)$ satisfy (LR), the dumbbell condition and metric doubling.
Does this space satisfy the EHI? 
}\end{problem}

\med {\bf Acknowledgment.} We are grateful to Laurent Saloff-Coste for some references and conversations 
related to this work, and Jun Kigami for some illuminating conversations on quasisymmetry.

\noindent MB: Department of Mathematics,
University of British Columbia,
Vancouver, BC V6T 1Z2, Canada. \\
barlow@math.ubc.ca

\noindent MM: Department of Mathematics, University of British Columbia and Pacific institute for the Mathematical Sciences,
Vancouver, BC V6T 1Z2, Canada. \\
mathav@math.ubc.ca

\end{document}